\documentclass[final]{siamltex}
\pdfoutput=1
\usepackage{graphicx,color,amsmath,amssymb,comment}
\usepackage{cite}
\usepackage[hypertexnames=false,colorlinks=true,linkcolor=blue,citecolor=blue,
            bookmarksopen=false,bookmarks=false,
            pdfstartview=XYZ,pdffitwindow=true,pdfcenterwindow=true]{hyperref}
\usepackage{braket}

\newcommand{\R}{\mathbb{R}}

\newcommand{\vect}[1]{\mathbf{#1}}
\newcommand{\matr}[1]{\mathbf{#1}}
\newcommand{\mr}{\vect{r}}

\newcommand{\dsix}{\mathbf{D6}}
\newcommand{\dfour}{\mathbf{D4}}
\newcommand{\dthree}{\mathbf{D3}}
\newcommand{\dtwo}{\mathbf{D2}}

\newcommand{\etal}{\textit{et al. }}

\newcommand{\frl}[1]{F^R_{l#1}}
\newcommand{\frr}[1]{F^R_{r#1}}
\newcommand{\fdsl}[1]{F^{\mathbf D6}_{l#1}}
\newcommand{\fdsr}[1]{F^{\mathbf D6}_{r#1}}
\newcommand{\pds}[0]{P^{\mathbf D6}}
\newcommand{\pdsfull}[0]{\pds_\textrm{IM}}

\newcommand\sref[1]{Sec.~\ref{#1}}

\newcommand\fref[1]{Fig.~\ref{#1}}

\newcommand\Fref[1]{Figure~\ref{#1}}
\newcommand\Frefs[1]{Figures~\ref{#1}}
\newcommand\eqreff[1]{Eq.~(\ref{#1})}
\newcommand\eqrefs[1]{Eqs.~(\ref{#1})}

\newtheorem{remark}[theorem]{Remark}
\title{Continuation of localised coherent structures in nonlocal
  neural field equations}

\author{
James Rankin\thanks{NeuroMathComp Team, INRIA Sophia Antipolis, France. E-mail: james.rankin@inria.fr} \and
Daniele Avitabile\thanks{Centre for Mathematical Medicine and Biology, School of
Mathematical Sciences, University of Nottingham, UK} \and 
Javier Baladron\thanks{NeuroMathComp Team, INRIA Sophia Antipolis, France} \and
Gregory Faye\thanks{School of Mathematics, University of Minnesota,  USA} \and 
David J.B. Lloyd\thanks{Department of Mathematics, University of Surrey, UK}}

\begin{document}

\maketitle

\begin{abstract}
  We study localised activity patterns in neural field equations posed
  on the Euclidean plane; such models are commonly used to describe
  the coarse-grained activity of large ensembles of cortical neurons
  in a spatially continuous way. We employ matrix-free Newton-Krylov
  solvers and perform numerical continuation of localised patterns
  directly on the integral form of the equation. This opens up the
  possibility to study systems whose synaptic kernel does not lead to
  an equivalent PDE formulation. We present a numerical bifurcation
  study of localised states and show that the proposed models support
  patterns of activity with varying spatial extent through the
  mechanism of homoclinic snaking. The regular organisation of these
  patterns is due to spatial interactions at a specific scale
  associated with the separation of excitation peaks in the chosen
  connectivity function. The results presented form a basis for the
  general study of localised cortical activity with inputs and, more
  specifically, for investigating the localised spread of orientation
  selective activity that has been observed in the primary visual
  cortex with local visual input.
\end{abstract}

\section{Introduction}

One of the most challenging research questions in neuroscience is
understanding the relationship between spatially-structured cortical
states and the underlying neural circuitry that supports them. A
popular approach for analysing coarse-grained activity of large
ensembles of neurons in the cortex is to model cortical space as a
continuum. Since the pioneering work of Wilson and Cowan
\cite{wilson-cowan:72,wilson-cowan:73} and Amari
\cite{amari:75,amari:77}, continuous neural field models have become a
popular and effective tool in neuroscience. In such models, the
large-scale activity of spatially-extended networks of neurons is
described in terms of nonlinear integro-differential equations, whose
associated integral kernels represent the spatial distribution of
neuronal synaptic connections. The canonical
Wilson-Cowan-Amari neural field equation
\cite{wilson-cowan:73,amari:77}
\begin{equation}\label{eq:wilsonCowan}
 \frac{\partial}{\partial t}
 u(\mr,t)=-u(\mr,t)+\int_{\Omega}w(\mr,\mr^\prime)S\big(u(\mr^\prime,t)\big)\text{dm}(\mr^\prime)
\end{equation}
describes the evolution of the average membrane voltage potential of a
neuronal population $u(\mr,t)$, at position $\mr\in\Omega$ on the
cortex and at time $t$. The nonlinear function $S$ represents the neural
firing rate, whereas the connectivity function $w(\mr,\mr')$ models
how a population of neurons at position $\mr$ on the cortex interacts
with a population at position $\mr^\prime$. Frequently-used firing
rate functions $S$ include the Heaviside step function,
piecewise-linear functions or smooth sigmoidal functions. Various
choices are also possible for the connectivity function, which is
often assumed to be translation invariant (that is, dependent on the
Euclidean distance $\Vert \mr - \mr^\prime \Vert$) and localised in
space. The cortical domain $\Omega$ is usually a subset of $\R^d$,
with $d=1$ or, in more realistic models, $d=2$. For a recent review on
neural fields modeling, we refer to Bressloff \cite{bressloff:12}.

Unlike spiking neural network models, continuous field models have the
advantage that analytic techniques for partial differential equations
(PDEs) can be adapted to study the formation of patterns and their
dependence upon control parameters. Various types of coherent
structures have been observed in neural field models, ranging from
spatially and temporally periodic patterns to travelling waves and
spiral waves~\cite{ermentrout:98,coombes:05,laing:13}. Neural field
equations have also successfully been used to model a wide range of
neurobiological phenomena such as visual hallucinations
\cite{ermentrout-cowan:79b,bressloff-cowan-etal:01}, mechanisms for
short term memory \cite{laing-troy-etal:02} and feature selectivity in
the visual cortex
\cite{ben-yishai-bar-or-etal:95,hansel-sompolinsky:97}.

A common strategy to derive analytical and numerical results for the
nonlocal~\eqreff{eq:wilsonCowan} is to assume translation invariance
and exploit the freedom in the choice of the connectivity function
$w$: if the Fourier Transform of $w$ is a rational function, it is
possible to derive a PDE formulation that is equivalent to the
integral model~\cite{laing:13}. Coherent structures supported by the
original model can then be conveniently constructed and analyzed in
the PDE framework. Indeed, previous studies have been carried out in
cases where the synaptic kernel led to an equivalent PDE
formulation~\cite{laing-troy-etal:02,laing-troy:03}.  To the best of
the author's knowledge, there has been no attempt to propose efficient
path-following methods for general connectivity functions $w$. This
paper is motivated by the desire to develop numerical algorithms
for~\eqreff{eq:wilsonCowan} without relying on an equivalent PDE
formulation.  More precisely, we discuss how to solve
\eqreff{eq:wilsonCowan} when $\Omega=\R^2$, $S$ is a smooth function
of sigmoidal type and $w$ has a generic Fourier Transform.

The main tools for our investigation are time simulation and numerical
continuation.  When a PDE model is available, we use standard
techniques for both tasks, in line with what is typically done in
several other works in this
field~\cite{laing-troy:03,laing:05,laing:13}. However we show that,
when the integral of \eqreff{eq:wilsonCowan} can be written as a
convolution, it is convenient to employ a fast Fourier transform (FFT)
for both time stepping and numerical continuation. Direct numerical
simulations of \eqreff{eq:wilsonCowan} using FFTs have been performed
before on full integral
models~\cite{folias-bressloff:05,coombes-schmidt-etal:12}. In the
present paper, we combine FFTs with Newton-Krylov
solvers~\cite{kelley:95,kelley:03}, thus opening up the possibility to
perform numerical continuation directly on the integral model. 

We concentrate our effort on the emergence and bifurcation structure
of stationary localised patterns in planar neural field models of the
form \eqreff{eq:wilsonCowan}. Indeed, this type of solution is of
great interest in models of prefrontal cortex, where localised states
are believed to characterise working memory
\cite{colby-duhamel-etal:95,funahashi-bruce-etal:89}.  Recently,
localised regions of activity have been observed in the cat primary
visual cortex \cite{chavane-sharon-etal:11} when the animal is
presented with localised-oriented input. Moreover, some reported
drug-induced visual hallucinations have also been found to be
spatially localised \cite{siegel:77} indicating the existence of
spatially localised regions of activity in the human primary visual
cortex.

Localised states have been observed in a wide variety of nonlinear
media~\cite{knobloch:08}. The bifurcation structure of localised
solutions has been studied extensively in the Swift--Hohenberg
equation posed on the real
line~\cite{coullet-riera-etal:00,woods-champneys:99,
  burke-knobloch:07,burke-knobloch:07c,chapman-kozyreff:09} and on the
plane~\cite{beck-knobloch-etal:09,lloyd-sandstede-etal:08,
  mccalla-sandstede:10,avitabile-lloyd-etal:10,sakaguchi-brand:98,
  sakaguchi-brand:96,lloyd-sandstede:09}.  In this context, a
well-known mechanism for the formation of localised states is
\textit{homoclinic snaking}: solutions with one or more bumps at the
core emerge from the trivial homogeneous state and undergo a series of
fold bifurcations, giving rise to a hierarchy of states with an
increasing number of bumps. This scenario seems to be a common
footprint of localised patterns, extending far beyond the prototypical
Swift--Hohenberg
equation~\cite{schneider-gibsol-etal:10,haudin-rojas-etal:11} and to
have also been found in nonlocal equations such as neural
fields~\cite{laing-troy-etal:02,coombes-lord-etal:03,
  laing-troy:03,faye-rankin-etal:12}.

\begin{figure}
\centering
\includegraphics{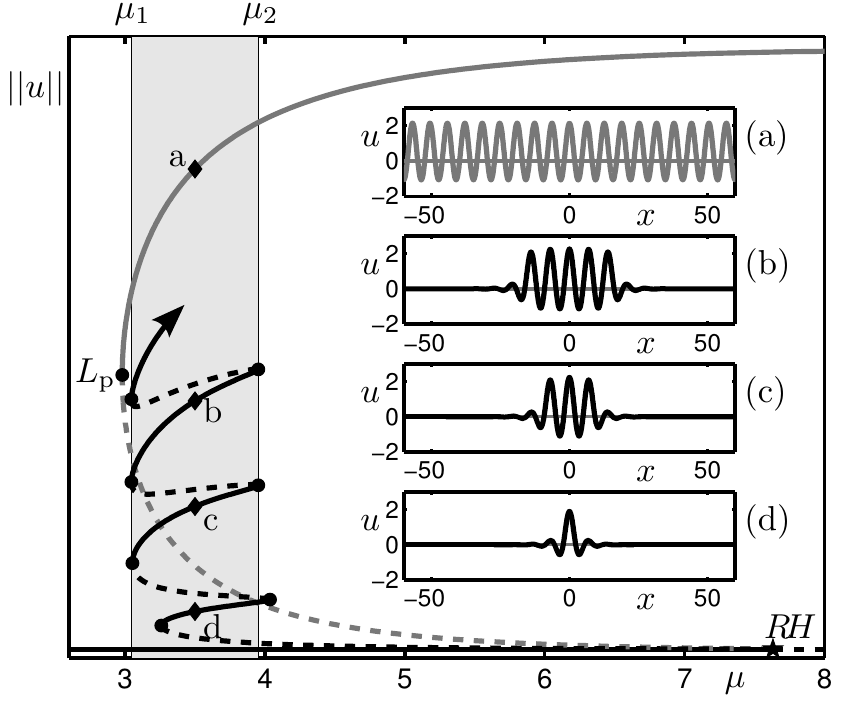}
\caption{Bifurcation diagram showing snaking behavior for the 1D
  neural field equation using~\eqrefs{eq:1D_w_fr}--(\ref{eq:1D_S_fr})
  with $\theta=3.5$ and $b=0.4$, where stable (unstable) branch
  segments are represented with solid (dashed) lines. Two branches of
  solutions bifurcate from the trivial states at a Reversible Hopf
  bifurcation $R\!H$, namely a branch of periodic solutions (grey) and
  a branch of localised solutions (black). The periodic branch is
  stable above the fold $L_{\textrm{p}}$. The localised branch
  undergoes a series of fold bifurcations giving rise to stable branch
  segments with increasing numbers of bumps. Several examples of
  solutions are shown in the insets (a)--(d). Stable localised states
  exist for $\mu\in[\mu_1,\mu_2]$.}
\label{fig:lt1d_snake}
\end{figure}

An example of homoclinic snaking is given in~\fref{fig:lt1d_snake},
where we show a bifurcation diagram for the integral
model~\eqref{eq:wilsonCowan} posed on the real line with
\begin{align}
  w(x,x'):=w(|x-x'|)&=e^{-b|x-x'|}(b \sin |x-x'| + \cos (x-x') ),\label{eq:1D_w_fr} \\
  S(u)&=\frac{1}{ 1+e^{-\mu u+\theta} }-\frac{1}{1+e^{\theta}}, \label{eq:1D_S_fr}
\end{align}
where $b$ and $\mu$ control the decay of the synaptic kernel and the
slope of the sigmoidal firing rate respectively, while $\theta$ is a
threshold value. As $\mu$ varies, the trivial steady state $u=0$
bifurcates at a subcritical Reversible Hopf bifurcation, from which a
branch of periodic solutions and a branch of localised states
originate. Beyond the fold $L_{\textrm{p}}$ on the periodic branch, a
stable periodic state, shown in the inset (a), coexists with the
trivial state. The branch of localised states features solutions with
an odd number of bumps and snakes for $\mu \in [\mu_1,\mu_2]$; in this
interval, localised solutions with different spatial extent coexist
and are stable (see insets (b)--(d)).  We note the existence of a
counterpart even-numbered-bump solution branch along with so-called
\textit{ladder branches} connecting the odd and even branches (not
shown). For the same connectivity function used here, snaking has been
shown to occur in terms of the parameter
$b$~\cite{laing-troy-etal:02}.  Elsewhere, Elvin \etal
\cite{elvin-laing-etal:10} used the Hamiltonian structure of the
steady states of the model \eqref{eq:wilsonCowan}
with~\eqref{eq:1D_w_fr}--\eqref{eq:1D_S_fr} and developed numerical
techniques to find homoclinic orbits of the system. Snaking has also
been shown to occur with a Mexican-hat connectivity
function~\cite{coombes-lord-etal:03} and for a wizard-hat connectivity
function~\cite{faye-rankin-etal:12}. In the latter study, normal form
theory for a Reversible Hopf bifurcation was applied to prove the
existence of localised solutions and a comprehensive parameter study
was carried out with numerical continuation in terms of two parameters
controlling the nonlinearity and a third controlling the shape of the
connectivity function.
\begin{figure}
  \centering
  \includegraphics[width=0.8\linewidth]{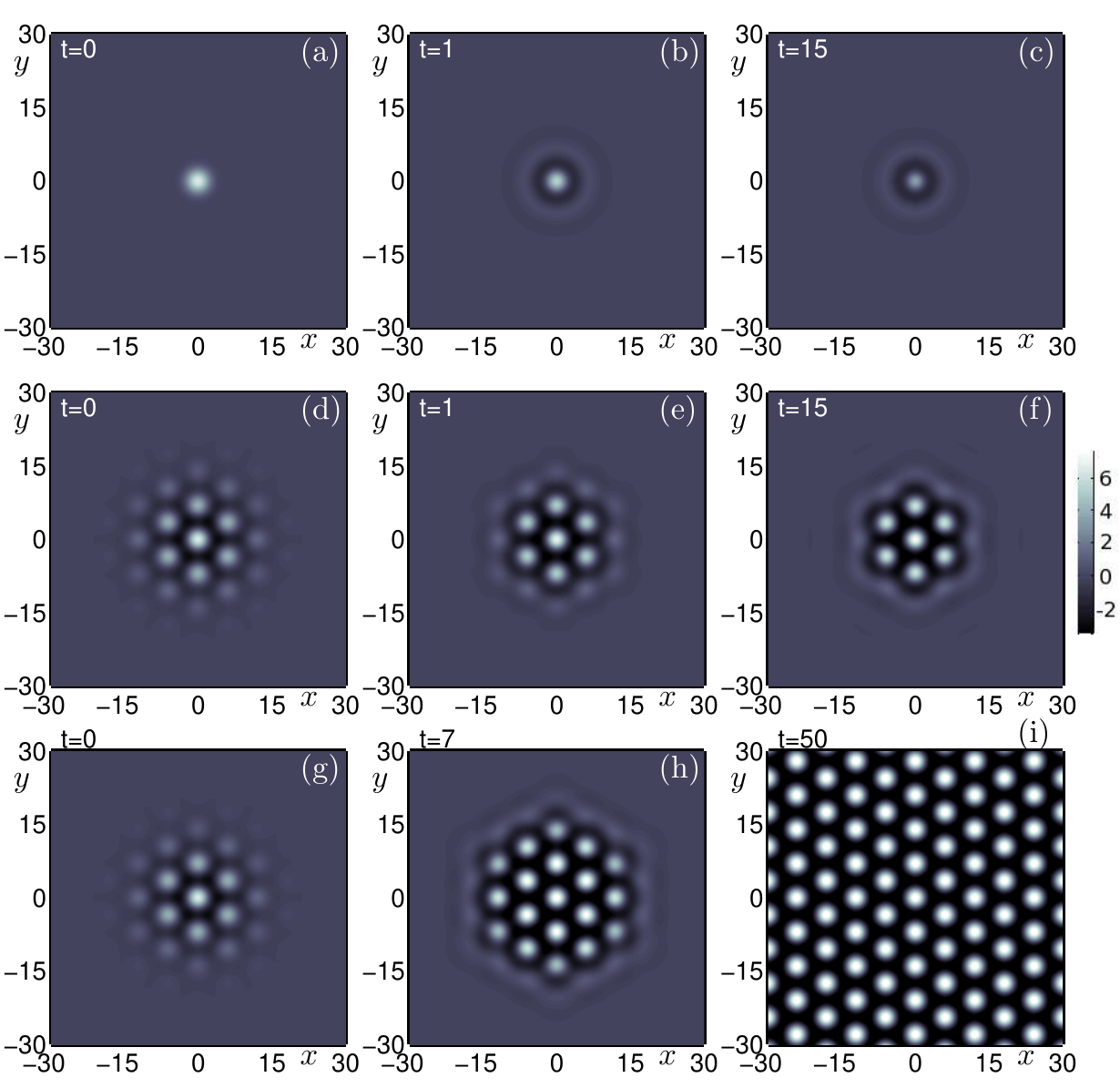}
  \caption{Time simulations of the integral
    model~\eqref{eq:wilsonCowan} posed on the plane, with connectivity
    function given by~\eqreff{eq:2D_w_fr} and firing rate
    function given by \eqreff{eq:1D_w_fr}; in all simulations
    $\theta=5.6$, $b=0.4$. (a)--(c): Convergence of a small bump of
    activity, given by \eqreff{eq:radic}, to a spot solution with
    $\mu=3.4$; time as indicated in panels. (d)--(f): Convergence of a
    hexagonal pattern, given by \eqreff{eq:hexic}, to a stable
    $\dsix$-symmetric localised state with $\mu=3.2$.  (g)--(i) With
    the same initial condition as (d), divergence away from a
    localised state to a periodic state with $\mu=3.4$.}
  \label{fig:full_th}
\end{figure}

For the neural field equation posed on the Euclidean plane, various
types of spatially localised two-dimensional states have been found including radially symmetric solutions
\cite{taylor:99,werner-richter:01,laing-troy-etal:02,
  laing-troy:03,folias-bressloff:05,bressloff-kilpatrick:11,
  faye-rankin-etal:12b},
rings~\cite{owen-laing-etal:07,coombes-schmidt-etal:12}, hexagonal
patches~\cite{laing-troy:03,faye-rankin-etal:12b} along with more
complex breathing and travelling
states~\cite{folias-bressloff:04,owen-laing-etal:07,
  coombes-schmidt-etal:13}. When working in the Euclidean plane, it is still possible to derive an
equivalent PDE for suitable choices of the connectivity function $w$
or of its Fourier transform $\widehat w$~\cite{laing-troy:03}.
Following this approach, normal form theory has recently been applied
to prove the existence of localised solution branches in both the
Euclidean plane and on the Hyperbolic disk with a wizard-hat
connectivity function~\cite{faye-rankin-etal:12b}.  When these
solutions were path-followed using numerical continuation, it was
found that the branches do not undergo snaking-type behaviour.
However, for the connectivity~\eqref{eq:1D_w_fr}, snaking was shown to
occur for branches of radially-symmetric
solutions~\cite{laing-troy:03}. Furthermore, the existence of
$\dsix$-symmetric and $\dthree$-symmetric localised states were found
at isolated parameter values.

In~\fref{fig:full_th} we show time simulations of the
nonlocal~\eqreff{eq:wilsonCowan} posed on the plane, with
radially-symmetric connectivity function
\begin{equation}\label{eq:2D_w_fr}
 w(\mr,\mr'):=w(\Vert \mr-\mr'\Vert)
             =e^{-b\Vert \mr-\mr'\Vert}(b \sin \Vert \mr-\mr'\Vert + \cos \Vert \mr-\mr'\Vert ),
\end{equation}
and sigmoidal firing rate function given by \eqreff{eq:1D_S_fr};
various combination of initial conditions and control parameters lead
to three different steady states. We note that this models supports
localised states, such as the radially-symmetric spot of panel (c) or
the hexagonal pattern of panel (f). Furthermore, changes in the slope
of the sigmoidal firing rate affect the stability properties of the
solutions, leading to other localised states or to domain-covering
patterns such as the one in panel (i).  In the present paper we will
focus on localised planar patterns (with various symmetry properties)
that coexist with the trivial state $u=0$ and with fully periodic
states, similar to the one shown in~\fref{fig:lt1d_snake} for the 1D
case.

A key result of the present paper is that homoclinic snaking occurs
in planar neural field models for non-radial patterns and that the
choice of the connectivity function has a considerable impact on the
snaking structure, as well as on the stability properties and the
selection of the localised states. In two spatial dimensions, periodic
and localised solution branches bifurcate from the trivial state at a
Turing instability (as opposed to a Reversible Hopf in 1D) and snake
irregularly, in a similar fashion to what is found for the planar
Swift--Hohenberg equation \cite{lloyd-sandstede-etal:08}.

The outline of the paper is as follows. In \sref{models}, we
present the different models that we study, and then, in
\sref{methods}, we describe the numerical methods that are used to
analyze each model. Our homoclinic snaking results of localised states
are presented in \sref{sec:numres}.

\section{Models}\label{models}
In this section, we introduce several neural field models used in the paper: our
starting points are the models introduced by Laing \etal
\cite{laing-troy:03,laing:13}, in which the cortical space $\Omega$ is assumed to be the
Euclidean plane $\R^2$. 

\subsection{Integral model}
\begin{figure}
\centering
\includegraphics[width=8cm]{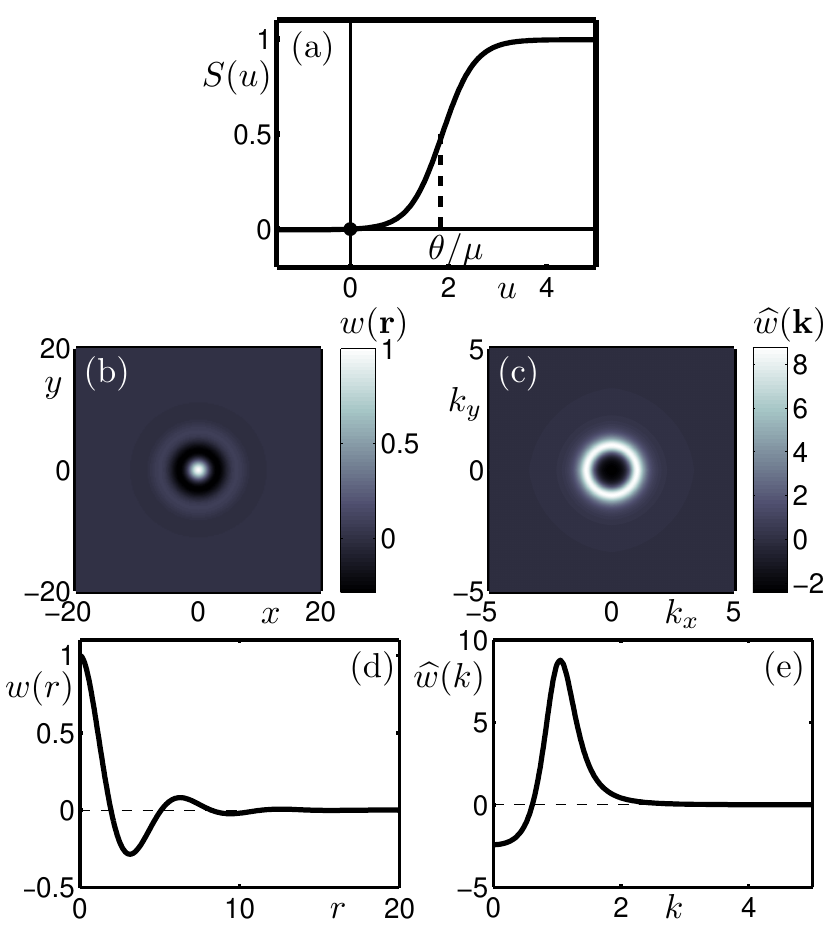}
\caption{Firing rate function and connectivity function for the
  integral model (IM). (a): Sigmoidal nonlinearity given by
  \eqreff{eq:intModelFR} plotted here with $\theta=5.6$ and $\mu=3$;
  in this formulation the effective threshold is $\theta/\mu$ (dashed
  line). (b): Radially-symmetric connectivity kernel given by
  \eqreff{eq:intModelW} plotted on the Euclidean plane with
  $b=0.4$. (c): Its Fourier transform plotted on the
  $(k_x,k_y)$-plane. (d): Radial cross section of panel (b). (e):
  Radial cross section of panel (c).}
\label{fig:full_conn_sig}
\end{figure}

The first model that we will consider is obtained from \eqreff{eq:wilsonCowan}
assuming a translation-invariant, radially-symmetric kernel
\begin{equation} \label{eq:intModel}
  \frac{\partial}{\partial t} u(\mr,t) = - u(\mr,t) + \int_{\R^2}
  w(\Vert \mr-\mr^\prime \Vert)\; S\big(u(\mr^\prime,t)\big)\; d\mr' + g(\mr)
\end{equation}
with sigmoidal firing rate (shown in \fref{fig:full_conn_sig}(a))
\begin{equation} \label{eq:intModelFR}
  S(u)=\frac{1}{ 1+e^{-\mu u+\theta} }-\frac{1}{1+e^{\theta}}, \qquad \mu,\theta >0,
\end{equation}
radial connectivity function
\begin{equation} \label{eq:intModelW}
  w\left( r \right)=e^{-b r} ( b \sin r+ \cos r ), \qquad r=\Vert \mr \Vert,
  \quad b > 0
\end{equation}
and external inhomogeneous input
\begin{equation}
  g(\mr) = G_0 \exp \Bigg( -\frac{\alpha x^2 + \beta y^2}{\sigma^2} \Bigg),
  \qquad G_0, \sigma \in \R, \quad \alpha, \beta > 0.
  \label{eq:extInput}
\end{equation}
In the visual cortex regions of activation have been shown to have a
Gaussian spread for radially-symmetric visual
inputs~\cite{chavane-sharon-etal:11}, hence our choice for the
function $g(\mr)$.
The connectivity function is plotted in \fref{fig:full_conn_sig}(b) on
the Euclidean plane and as a radial cross section
in~\fref{fig:full_conn_sig}(d). An analytical expression for the
Fourier transform of $w$ cannot be obtained; we show the Fourier
transform of \eqreff{eq:intModelW} as computed numerically in the
$\mathbf{k}=(k_x,k_y)$-plane in~\fref{fig:full_conn_sig}(c) and as a
radial cross section in~\fref{fig:full_conn_sig}(e), where $k=\Vert
\mathbf{k}\Vert$. The connectivity function~\eqref{eq:intModelW} was
proposed in models of working memory as a description of synaptic
connections in the prefrontal
cortex~\cite{gutkin-ermentrout-etal:00,laing-troy-etal:02}. The
connectivity describes local excitation and longer-range connections
that alternate between inhibition ($w(r)<0$) and excitation
($w(r)>0$). We argue that this type of connectivity pattern is also
relevant to the study of patterns of activity in early visual areas
like V1 where there is a characteristic length scale associated with
the average orientation hypercolumn width.  It has been shown in
anatomical studies that the number of lateral connections decay with
distance, that the number of excitatory connections peak each
hypercolumn width and the number of inhibitory connections peak each
half-hypercolumn width~\cite{buzas-eysel-etal:01}.  The net effect is
alternating bands of inhibition and excitation that decay with
distance. This is also consistent with auto-correlations computed for
the orientation selectivity map~\cite{niebur-wortgotter:94} given that
connections tend to be reinforced between neurons with similar
orientation preference. Henceforth the
model~\eqref{eq:intModel}--\eqref{eq:intModelW} will be referred to as
the integral model (IM).

\subsection{Fourth-order PDE approximation} \label{sec:PDE4} In the
cases when the Fourier transform of the synaptic kernel is a rational
function, it is possible to derive an equivalent PDE formulation of
\eqreff{eq:intModel}~\cite{laing:13}. For simplicity, we will consider
models without an external input $g(\mr)$. If $\widehat w(k) =
P(k)/Q(k)$ with $P$ and $Q$ even functions in $k$ where $k=\Vert
\mathbf{k} \Vert $ for $\mathbf{k}\in\R^2$, then the Fourier transform
of \eqreff{eq:intModel} gives
\[
  Q(k) \bigg[ \frac{\partial}{\partial t} \widehat u(\mathbf{k},t) + \widehat u(\mathbf{k},t)
  \bigg] = P(k) \widehat{(S \circ u)}(\mathbf{k},t).
\]
An inverse Fourier transform of the equation above leads to the desired PDE
\[
 \mathcal{L}_Q \bigg[ \frac{\partial}{\partial t} u(\mr,t) + u(\mr,t) \bigg] = 
 \mathcal{L}_P S \big[ u(\mr,t) \big],
\]
where $\mathcal{L}_P$ and $\mathcal{L}_Q$ are linear operators containing spatial
derivatives of even order. 

Since the Fourier transform of the connectivity function~\eqref{eq:intModelW} does
not have an analytic expression, the integral model does not admit an equivalent PDE.
However, we can approximate $w$ with a function whose Fourier transform is rational
and then derive an \textit{approximate} PDE for the integral model. We specify the
approximate connectivity function $w_4(r)$ through its Hankel Transform
\[
w_4(r)=\frac{1}{2\pi}\int_0^\infty s\, \widehat{w_4}(s) J_0(rs) \;ds,
\]
where $\widehat{w_4}(k) \approx \widehat w(k)$ is given by
\begin{equation}\label{eq:what_2d_4th}
  \widehat{w_4}(k)=\frac{A}{B+(k^2-M)^2}
\end{equation}
and the coefficients $A$, $B$, $M$ are determined using a
least-squares best fit algorithm (see \sref{sec:lsqfit} for further
details).

\begin{figure}
\centering
\includegraphics{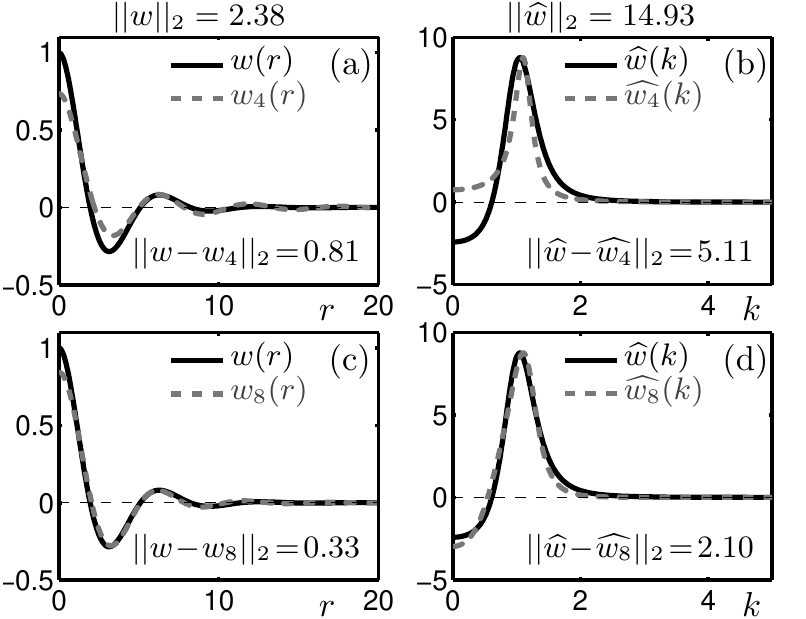}
\caption{Comparison of the PDE connectivity functions with $w(r)$
  given by~\eqreff{eq:intModelW} and its Fourier transform
  $\widehat{w}(r)$ as computed numerically. Panels (a) and (b) show
  the $4$th-order approximation defined
  by~\eqreff{eq:what_2d_4th} in the real and Fourier domains,
  respectively. Similarly, for the $8$th-order approximation 
  given by \eqreff{eq:what_2d_8th} in panels (c) and (d). For
  reference, we indicate the $L_2$-norm of $w$, $w-w_4$, $w-w_8$ and
  of their Fourier transforms.  Parameters given at the beginning of
  \sref{snaking1}.}
\label{fig:ltpde_conn}
\end{figure}
We compare the approximation with the original connectivity function
in the physical and Fourier domains in \fref{fig:ltpde_conn}(a) and
(b). In physical space the two functions appear to be similar. The key
qualitative difference is that in Fourier space $\widehat{w_4}(0)>0$,
which is not consistent with the IM connectivity function, for which
$\widehat{w}(0)<0$. Biologically this means that $w_4$ represents a
globally excitatory connectivity function, whereas $w$ represents a
globally inhibitory connectivity function. We will see in
\sref{sec:8th} that it is necessary to increase the order of the
polynomials in the numerator and denominator
of~\eqreff{eq:what_2d_4th} in order to accurately capture the sign at
$k=0$.

Starting from the expression for $\widehat{w_4}$, we derive the corresponding
PDE, containing spatial derivatives up to the fourth order
\begin{equation}\label{eq:PDE4}
  \big[ B+(M + \Delta)^2 \big] \bigg[ \frac{\partial}{\partial t} u(\mr,t) +
  u(\mr,t) \bigg] = A S \big( u(\mr,t) \big),
\end{equation}
where the sigmoidal firing rate function is identical to the integral
model case \eqref{eq:intModelFR}. In \eqreff{eq:PDE4} we have denoted
by $\Delta$ the standard Euclidean Laplacian.  Henceforth, this model
will be referred to as PDE4.

\subsection{Eight-order PDE approximation}
In order to obtain a more accurate representation of the connectivity
function~\eqref{eq:intModelW}, we repeat the steps outlined in
\sref{sec:PDE4} with the following approximation
\begin{equation}\label{eq:what_2d_8th}
\widehat{w_8}(k)=\frac{-A(k^2-C)(k^2-D)}{B+(k^2-M)^4},
\end{equation}
where the values of $A$, $B$, $C$, $D$ and $M$ are determined using a
least-squares best fit algorithm; see \sref{sec:lsqfit} for further
details. We compare this approximation with the original connectivity
function in the physical and Fourier domains in
\fref{fig:ltpde_conn}(c) and (d). With the higher-order polynomials
used here the approximation is more accurate and $\widehat{w_8}(0)>0$,
consistent with IM (compare the $8$th-order approximation as shown
in~\fref{fig:ltpde_conn}(c) and (d) with the $4$th-order approximation
as shown in~\fref{fig:ltpde_conn}(a) and (b)).  The synaptic function
$w_8$ leads to the following PDE, containing spatial derivatives up to
the eight order
\begin{equation}\label{eq:PDE8}
  \big[ B+(M + \Delta)^4 \big] \bigg[ \frac{\partial}{\partial t} u(\mr,t) +
  u(\mr,t) \bigg] = -\big[ (\Delta + C) (\Delta + D) \big] S \big( u(\mr,t)
  \big), 
\end{equation}
where the sigmoidal firing rate function is identical to the integral model
case. Henceforth this model will be denoted as PDE8.

\section{Numerical methods}\label{methods}
In this section, we review the numerical methods employed for the computation of
localised states in IM, PDE4 and PDE8.

\subsection{Integral model} Numerical computations of the
IM~\eqref{eq:intModel}--\eqref{eq:intModelW} are performed
discretizing a large but finite domain $\Omega = [-L,L]^2$ with $N$
evenly-distributed grid points in each spatial direction and imposing
periodic boundary conditions. We approximate $u$ on a grid $\Omega_N =
\set{(x_i,y_j)}_{i,j=1}^N$ and collect the corresponding approximate
values of $u$ in a vector $\vect{u}$
\[
  u_{ij} \approx u(x_i,y_j), \qquad \vect{u} = \set{u_{ij}}_{i,j=1}^N \in \R^{N^2}. 
\]
Similarly, we form vectors $\vect{w},\, \vect{S}(\vect{u}),\, \vect{g}
\in \R^{N^2}$ for the approximations to $w$, $S(u)$ and $g$,
respectively. Further, we introduce the discrete
convolution, 
\begin{equation}
  (u \ast v)_{ij} \approx \mathcal{F}^{-1} \big( \mathcal{F}(u) \,
  \mathcal{F}(v) \big) (x_i,y_j), \qquad
  \vect{u} \ast \vect{v} = \set{(u \ast v)_{ij}}_{i,j=1}^N \in \R^{N^2},
\label{eq:discrConv}
\end{equation}
where we have denoted by $\mathcal{F}$ and $\mathcal{F}^{-1}$ the 2D
Fourier Transform and its inverse, respectively.  In summary, the
discrete version of the evolution equation~\eqref{eq:intModel} is
given by
\begin{equation}
  \dot{\vect{u}} = - \vect{u} + \vect{w} \ast \vect{S}(\vect{u}) + \vect{g}.
  \label{eq:IMTimeStep}
\end{equation}
This type of discretization has been applied before in direct
numerical simulations of neural models (see, for
instance~\cite{folias-bressloff:05,coombes-schmidt-etal:12}) even
though it has not been used for numerical continuation. For smooth
firing rate functions, the right-hand side can be evaluated accurately
and efficiently using a Fast Fourier Transform (FFT) and its inverse
(IFFT). In passing, we note that since the FFT of $\vect{w}$ can be
performed and stored at the beginning of the computation, one function
evaluation of the right-hand side requires just one FFT and one
IFFT. Furthermore, standard de-aliasing techniques can be applied to
the convolution operator if required~\cite{canuto-hussaini:06}.

Once a stable steady-state of \eqreff{eq:IMTimeStep} is found via
direct numerical simulation, it is possible to continue it in one of
the control parameters using standard numerical continuation
techniques. In previous studies of neural field equations, numerical
continuation was performed on an equivalent PDE formulation of the
integral system. A key observation is that path following can be
applied directly to IM (or to similar models), employing FFTs and
Newton-Krylov solvers~\cite{kelley:95,kelley:03}. Such methods do not
require the formation of a Jacobian matrix, but rely only on
Jacobian-vector multiplications: for IM, this is conveniently done
using just a single application of FFT and IFFT. We remark that
Newton-Krylov methods are often used in conjunction with sparse
systems, but the performance of FFTs and IFFTs makes them a favourable
choice for IM, even though the system is full.

For numerical continuation of steady states of IM, we solve the system of
algebraic equations
\begin{equation}\label{eq:IMCont}
  \vect{F}(\vect{u}) =  - \vect{u} + \vect{w} \ast \vect{S}(\vect{u}) + \vect{g} = \vect{0},
\end{equation}
whose associated Jacobian-vector product is given by
\begin{equation}\label{eq:IMJac}
  \vect{J}(\vect{u}) \vect{v} =  - \vect{v} + \vect{w} \ast \big(
  \matr{S}^\prime(\vect{u}) \vect{v} \big),\quad
  \vect{u}, \vect{v} \in \R^{N^2},
\end{equation}
where $\matr{S}^\prime(\vect{u}) = \diag(S'(u_{11}),\ldots,S'(u_{NN}))
\in \R^{N^2 \times N^2}$.  We solve the system \eqref{eq:IMCont} using
a Newton-GMRES solver implemented in Matlab and continue the solution
with a secant method. Eigenvalue computations can also be performed
using the Jacobian-vector products \eqref{eq:IMJac}. Details of the
numerical implementation and of the numerical parameters can be found
in~\sref{subsubsec:IMImplementation}.

\begin{remark}
  The external input $g$ guarantees that the system of algebraic
  equations~\eqref{eq:IMCont} is not translation invariant, even when
  the problem is complemented with periodic boundary
  conditions. Unless otherwise stated, we will use a negligible
  external input for the IM, so that Newton iterations can be applied
  directly to~\eqreff{eq:IMCont}. We point out that a similar result
  can be obtained without imposing any external input, by perturbing
  the right-hand side with a term containing one extra unknown and
  then closing the system with a suitable phase
  condition~\cite{champneys-sandstede07,lloyd-sandstede-etal:08,avitabile-lloyd-etal:10}.
\end{remark}

\subsection{Fourth-order PDE approximation}
In order to continue stationary localised solutions to PDE4 with $\dsix$ symmetry, we
use polar coordinates and pose a boundary-value problem on the sector $\Omega_{1/6} =
\set{ (r,\theta) \in \R^2 | 0 < r < R,\; 0 < \theta < \pi/3}$ with Neumann boundary
conditions
\begin{equation}
\begin{aligned}
0 = &  \big[ B+(M + \Delta)^2 \big] u - A S(u),  && \quad (r,\theta) \in \Omega_{1/6} \\
0 = & \nabla u  \cdot \vect{n} && \quad (r,\theta) \in \partial \Omega_{1/6}
\end{aligned}
\label{eq:PDE4BVP}
\end{equation}
where the Laplacian operator $\Delta$ is expressed in polar coordinates. We
discretise the system above using finite differences in $r$ and a Fourier
collocation method in $\theta$, with $N_r$ and $N_\theta$ evenly spaced points,
respectively, leading to a system of nonlinear algebraic equations of the form
\begin{equation}\label{eq:PDE4Cont}
  \vect{0} = \big[ B\matr{I}+(M\matr{I} + \matr{L})^2 \big] \vect{u} - A \vect{S}(\vect{u})
  \qquad \vect{u} \in \R^{N_r N_\theta}, \quad \matr{I},\matr{L} \in \R^{N_rN_\theta
  \times N_rN_\theta}
\end{equation}
where $\matr{I}$ is the identity matrix. The Laplacian matrix $\matr{L}$ is formed
explicitly, starting from differentiation matrices $\matr{D}_r$, $\matr{D}_{rr}$,
$\matr{D}_{\theta\theta}$ for spatial derivatives with Neumann boundary conditions
and then combining them with Kronecker
products~\cite{trefethen:00,avitabile-lloyd-etal:10,lloyd-sandstede-etal:08}
\[
\matr{L} = \matr{D}_{rr} \otimes \matr{I_\theta} + ( \matr{R}^{-1} \matr{D}_r )
\otimes \matr{I_\theta} + \matr{R}^{-2} \otimes \matr{D}_{\theta\theta},
\]
where $\matr{R} = \diag(1,r_2,\ldots,r_{N_r})$ and $\matr{I}_\theta$ is the
$N_\theta$-by-$N_\theta$ identity matrix. For purely radial patterns, we adapt
the boundary-value problem so as to contain only the radial direction $r$ and
impose Neumann boundary conditions. Numerical continuation of the
system~\eqref{eq:PDE4Cont} is performed with a secant method. Further details on
the numerical implementation can be found
in~\sref{subsubsec:PDEnumericalDetails}.

\subsection{Eight-order PDE approximation}
For localised solutions of PDE8, we follow a similar approach to the
on used for PDE4. In order to avoid the discretisation of 8th-order
differential operators, we recast PDE8 as a system of two 4th-order
PDEs and seek solutions to the following boundary-value problem
\begin{equation}
\begin{aligned}
0 = &  (M + \Delta)^2 u - v,  && \quad (r,\theta) \in \Omega_{1/6}, \\
0 = &  (M + \Delta)^2 v + Bu + A (\Delta + C)(\Delta + D) S(u) ,  
&& \quad (r,\theta) \in \Omega_{1/6}, \\
0 = & \nabla u  \cdot \vect{n}, && \quad (r,\theta) \in \partial \Omega_{1/6}, \\
0 = & \nabla v  \cdot \vect{n}, && \quad (r,\theta) \in \partial \Omega_{1/6}.
\end{aligned}
\label{eq:PDE8BVP}
\end{equation}
Again, the discretisation of this system uses finite-differences in
$r$ and Fourier spectral collocation in $\theta$. The example above
shows also that, as we approximate $w$ more accurately, the order of
the underlying PDE increases and its numerical continuation becomes
more demanding. A more convenient approach would be to use directly
the integral form for the model, with connectivity function $w_8$ and
proceed with a Newton-GMRES, solver. In this way, the computational
cost would be the same as for IM.

\subsection{Implementation and numerical parameters} 

\subsubsection{Integral model}
\label{subsubsec:IMImplementation}
Time simulations are carried out using a standard $4$th order
Runge--Kutta method with fixed step size. At each time step, the
right-hand side of~\eqreff{eq:IMTimeStep} is evaluated four times
using an Nvidia Graphic Processing Unit (Tesla C2070). To compute the
Discrete Fourier Transform we use the CUFFT library provided by
\textsc{Nvidia} as part of its CUDA framework
\cite{nickolls-buck-etal:08}. This software library allows us to
easily exploit the parallelism of a GPU to obtain a fast
implementation. The vector $\vect{u}$ is kept in the GPU global memory
throughout a time step in order to avoid memory transfers and it is
updated in parallel once the four stages of the Runge--Kutta scheme
have been computed. Transfers between CPU and GPU memory only occur
when the result of a time step needs to be saved to a file. Time
simulations use a grid of $10^3 \times 10^3$ points and a stepsize of
$0.5$. Computation of each time step requires approximately $0.1$
seconds.

Numerical continuation for IM is done in Matlab using a in-house
secant continuation code which employs a Newton-GMRES method for the
nonlinear solves. Unless otherwise stated, we used a grid of $2^{10}
\times 2^{10}$ points and fixed an absolute tolerance of $10^{-3}$ for
the nonlinear iterations. The Newton-GMRES solver uses the MATLAB
in-built function \texttt{gmres} without preconditioners and with
parameters $\texttt{restart} = 20$, $\texttt{tol} = 10^{-3}$ and
$\texttt{maxit} = 10$. Initial guesses are obtained directly from the
Runge--Kutta time stepper, interpolating with MATLAB's
\texttt{interp2} function where necessary. Stability computations are
performed with Arnoldi iterations via MATLAB's \texttt{eigs} function,
passing the Jacobian-vector product~\eqref{eq:IMJac} and computing
(with the default tolerance) the first $20$ eigenvalues with the
largest real part. Computations are performed on a MacPro with a $3$
GHz Quad-Core Intel Xeon processor employing exclusively the CPU.

\subsubsection{PDE4 and PDE8 models}
\label{subsubsec:PDEnumericalDetails} Numerical continuation for the
PDE models have been carried out with a secant code similar to the one
used for IM, but using MATLAB's in-built function \texttt{fsolve} for
the nonlinear iterations. Unless otherwise stated, we used $300$ grid
points in the radial directions and $20$ in the angular direction. We
use the Levenberg--Marquardt algorithm implemented in \texttt{fsolve}
and set $\texttt{TolFun} = 10^{-6}$. The sparse Jacobians of these
problems are formed and passed directly to the solver. Initial guesses
for the continuation have been obtained using the expressions given in
\eqrefs{eq:radic} and \eqref{eq:hexic}. Computations are performed on
a standard laptop on a single core.

\subsubsection{Least-squares data fitting}\label{sec:lsqfit}
Before comparing the connectivity functions $w_4(r)$ and $w_8(r)$ with
$w(r)$, it was necessary to tune the parameters in the definitions of
$\widehat{w_4}$ and $\widehat{w_8}$. In each case we performed a
nonlinear least-squares optimization of the parameters using the
\texttt{lsqcurvefit} function in Matlab. For $w_4$ the objective was
to minimize the $L_2$-norm of the difference between $\widehat{w_4}$
and $\widehat{w}$ whilst varying the parameters $A$, $B$ and $M$ in
\eqreff{eq:what_2d_4th}, where $\widehat{w}$ is computed numerically
using the Hankel transform at $300$ points.  Similarly, for $w_8$, the
$L_2$-norm of the difference between $\widehat{w_8}$ and $\widehat{w}$
was minimised whilst varying the paramaters $A$, $B$, $M$, $C$, $D$
in~\eqreff{eq:what_2d_8th}.  For reference, the $L_2$-norms of $w$ and
$\widehat{w}$ are given above the panels in \fref{fig:ltpde_conn}; the
norm of the difference between the two functions plotted in each panel
is also given. We note the largest Fourier mode of the connectivity
dictates the location of bifurcations in terms of the sigmoid
parameters $\theta$ and $\mu$. By minimising the difference between
the connectivity functions in Fourier space, we expect to find similar
behaviour for each connectivity over the same parameter ranges. On the
other hand, if one minimises the difference between the connectivity
functions in physical space and the amplitudes of the largest Fourier
modes are not matched, bifurcations occur in different parameter
ranges in each model and a direct comparison cannot be made.

\section{Numerical results}\label{sec:numres}
\subsection{Convergence of the Newton-GMRES solver}\label{gmres}
\begin{figure}
\centering
\includegraphics{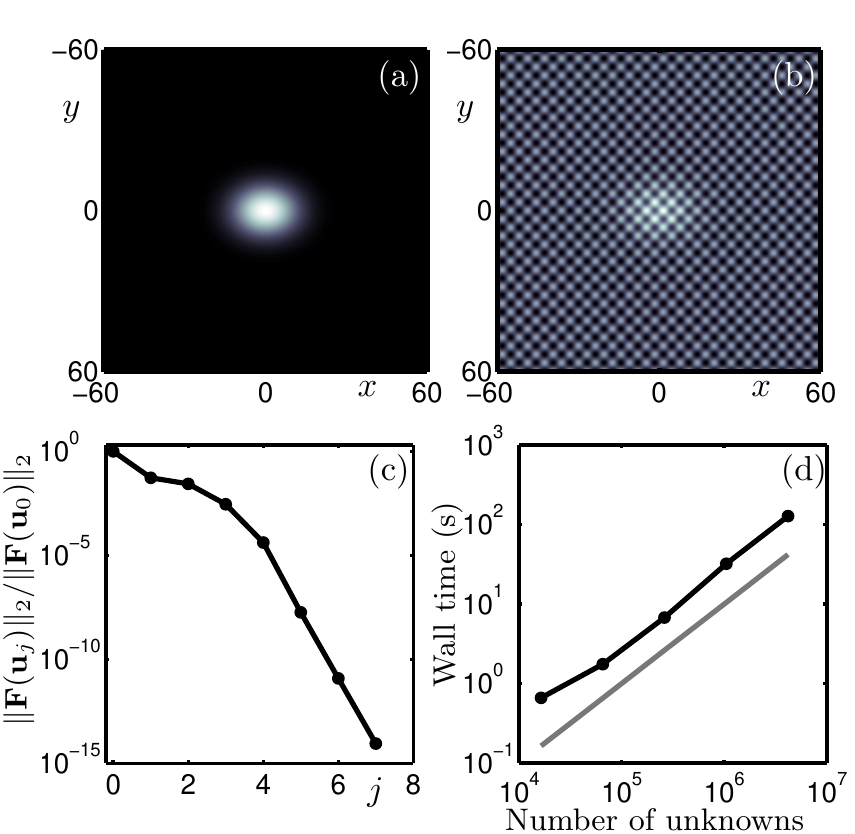}
\caption{Convergence of the Newton-GMRES solver. A stationary
  localised solution $u_\ast$ of IM is shown in panel (a). The
  solution is perturbed to obtain an initial guess $u_0 = u_\ast + 0.8
  \sin x \cos y$, shown in panel (b), and then the Newton-GMRES solver
  is used to converge back to $u_\ast$. Convergence is measured with
  the relative residual $\Vert \vect{F}(\vect{u}_j) \Vert_2 / \Vert
  \vect{F}(\vect{u}_0) \Vert_2$, where $\vect{u}_j$ is the solution at
  the $j$th iteration (panel (c)). The experiment is repeated for
  several values of $N$ without any significant change to the
  convergence diagram, whereas the wall time for the numerical experiment
scales linearly with the number of unknowns $N^2$, as reported in panel (d).
Parameters of IM: 
$\theta  = 5.6  $, 
$\mu     = 2.5  $, 
$b       = 0.40 $,
$L       = 60   $,
$G_0      = 4.0 $,
$\alpha  = 1.0  $,
$\beta   = 4.0  $,
$\sigma  = 12.0 $.}
\label{fig:newtongmres}
\end{figure}

Since Newton-GMRES methods with pseudospectral evaluation of the
right-hand side have not been used before for integral neural field
models, we report briefly on our solver. To test convergence, we
perturbed a localised steady state of IM to obtain an initial guess
(panel (b) of \fref{fig:newtongmres}) and converge back to the
original solution (panel (a) of \fref{fig:newtongmres}) using our
Newton-GMRES solver. In panel (c) we plot the relative residuals of
each iteration, showing that we achieve convergence within a few
nonlinear iterations. Similar convergence plots (not shown) are
obtained for the numerical continuation, albeit solutions in that case
are achieved with fewer iterations, owing to the more accurate initial
guess provided by the secant predictor scheme.  The experiment is
repeated for various values of $N$: the convergence diagrams are
indistinguishable from the one reported in panel (c). The wall time
for the numerical experiment scales linearly with the number of
unknowns $N^2$, as reported in panel (d). We remark that, even without
using any GPU acceleration and without enforcing explicit
parallelisation in the CPU, the Newton-GMRES solver finds a solution
to a full problem with $1,\!048,\!576$ unknowns in less than $40$ seconds.

\subsection{Snaking behaviour of radial and $\dsix$
  patterns}\label{snaking1}
We now turn to the numerical continuation of localised states in IM,
PDE4 and PDE8. The continuation parameter is the steepness $\mu$ of
the sigmoidal firing rate, whereas the other parameters are fixed as
follows: for IM, we choose $\theta = 5.6 $, $b = 0.40 $, $L = 60 $,
$G_0 = 10^{-4} $, $\alpha = 1.0 $, $\beta = 1.0 $, $\sigma = \sqrt{10}
$; for the PDE4 model, we use $\theta = 5.6$, $R = 60$, $G_0 = 0$ with
fitting parameters for $\widehat{w_4}$ in~\eqreff{eq:what_2d_4th}
given by $A=1.225$, $B=0.1398$, $M=1.2183$; for the PDE8 model we use
again $\theta = 5.6$, $R = 60$, $G_0 = 0$ with fitting parameters for
$\widehat{w_8}$ in~\eqreff{eq:what_2d_8th} given by $A=0.8510$,
$B=0.6626$, $M = 0.6653$, $C=0.3$ and $D=10$. We remark that
translation invariance is removed in the IM by the negligible external
input $G_0$ while in PDE4 and PDE8 this is achieved by the boundary
conditions of the problems~\eqref{eq:PDE4BVP} and~\eqref{eq:PDE8BVP},
so we choose $G_0=0$. In the present section we focus on the no- (or
equivalent negligible-) input case which should be well understood
before the addition of an input. In~\sref{sec:discuss_inp} we will
provide examples of the model behaviour with input and discuss the
implications.

The bifurcation points in the diagrams presented in this section will
be labelled as follows. $F$ represents a fold bifurcation and $P$ a
spatial-symmetry-breaking bifurcation from a radial state.
Superscripts indicate the symmetry properties of the bifurcation,
where $R$ represents a bifurcation on a branch with radially symmetric
solutions and $\dsix$ represents a bifurcation on a branch of
$\dsix$-symmetric solutions The labels $l$ and $r$ in the subscripts
for fold bifurcations indicate whether the fold occurs on the left or
right of the snaking structure. The indices $n$ in the subscripts
indicate the ordering moving up the snaking structure. On radial
branches $n$ corresponds to the number of rings around a central spot
solution, for example, on the branch between $\frl{1}$ and $\frr{1}$
there is one ring around a central spot. On $\dsix$-branches there are
$n(n+1)/2\times6$ additional spots glued around a central spot, for
example, on the branch between $\fdsl{1}$ and $\fdsr{1}$ the solution
has a total of 7 spots.

\subsubsection{PDE4 results and comparison with IM}\label{sec:4th}

\begin{figure}
\centering
\includegraphics[width=.8\linewidth]{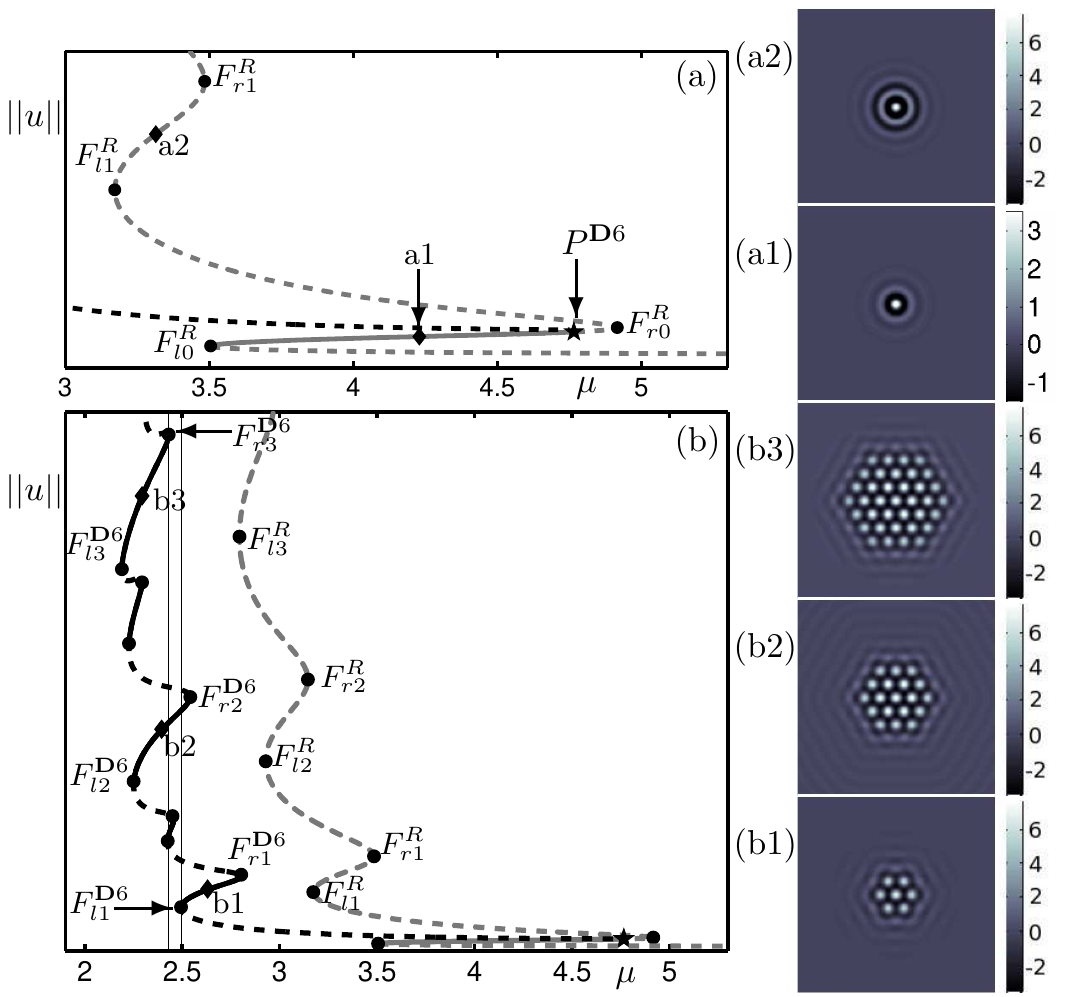}
\caption{Snaking of radially- and $\dsix$-symmetric localised
  solutions in terms of $\mu$ for PDE4. Radial branches are grey
  curves and $\dsix$ branches are black; stable segments are solid and
  unstable segments are dashed. Bifurcations labelled $F$ and $P$ are
  discussed in the text. (a): Detail of radial branch from (b); points
  labelled a1 and a2 correspond to planar plots showing a stable spot
  solution and an unstable spot-with-ring solution, respectively. (b):
  Global structure showing radial and $\dsix$ branches; points
  labelled b1, b2 and b3 correspond to planar plots showing $7$-spot,
  $19$-spot and $37$-spot solutions, respectively.  Thin vertical
  lines discussed in text.  Parameters given at the beginning of
  \sref{snaking1}.}
\label{fig:lt4th_snake_1} 
\end{figure}

We discuss both radially- and $\dsix$-symmetric localised solutions of
PDE4 with connectivity function defined via
\eqreff{eq:what_2d_4th}. Other solution branches with different
symmetry properties do exist but these are only discussed for IM in
\sref{sec:numerics_full}. An unstable radial spot solution bifurcates
from the trivial state $u=0$ at a Turing instability with $\mu$-value
to the right of the range shown in the subsequent diagrams. It is this
unstable radial spot branch that appears in the bottom right of
\fref{fig:lt4th_snake_1}(a) and (b).

\Fref{fig:lt4th_snake_1}(b) shows the snaking structure for radial and
$\dsix$ branches. We first focus on the radial branch, detail of which
is shown in panel (a). A radial spot branch enters the diagram in the
bottom-right-hand corner and undergoes a fold at $\frl{0}$. The radial
spot solution existing on the branch segment between $\frl{0}$ and
$\frr{0}$ is plotted in panel (a1). This solution is stable between
$\frl{0}$ and $\pds$. At $\pds$ on the radial branch a $\dsix$
instability produces the bifurcating branch of $\dsix$-symmetric
solutions that leaves panel (a) in the bottom-left-hand corner. Beyond
$\pds$ the radial branch is unstable and undergoes a further fold
$\frr{0}$. After another fold $\frl{1}$ a ring has formed around the
radial spot. The spot with ring solution existing on the branch
segment between $\frl{1}$ and $\frr{1}$ is shown in panel (a2). The
branch remains unstable and undergoes a series of further folds
($\frl{2}$, $\frr{2}$, $\frl{3}$, etc) adding additional rings as is
shown in panel (b).

The $\dsix$-symmetric branch that bifurcates from the radial branch at
$\pds$ also undergoes a series of fold bifurcations as shown in
\fref{fig:lt4th_snake_1}(b). In this case a series of additional spots
are added to the central spot in a configuration that preserves the
$\dsix$-symmetry. Planar plots in panels (b1), (b2) and (b3) show the
stable $7$-spot, $19$-spot and $37$-spot solutions that exist on the
branch segments between the fold pairs $(\fdsl{1},\fdsr{1})$,
$(\fdsl{2},\fdsr{2})$ and $(\fdsl{3},\fdsr{3})$, respectively. There
are further intermediate stable branch segments between pairs of fold
bifurcations that have not been labelled. On the stable segment that
can be found between $\fdsr{1}$ and $\fdsl{2}$, there exists a
$13$-spot solution for which one spot has been glued on the long edge
of each of the six sides of the solution shown in panel
(b1). Similarly, on the stable segment that
can be found between $\fdsr{2}$ and $\fdsl{3}$, there exists a
$25$-spot solution for which two spots have been glued on the long edge
of each of the six sides of the solution shown in panel (b2).

\begin{figure}
\centering
\includegraphics[width=8cm]{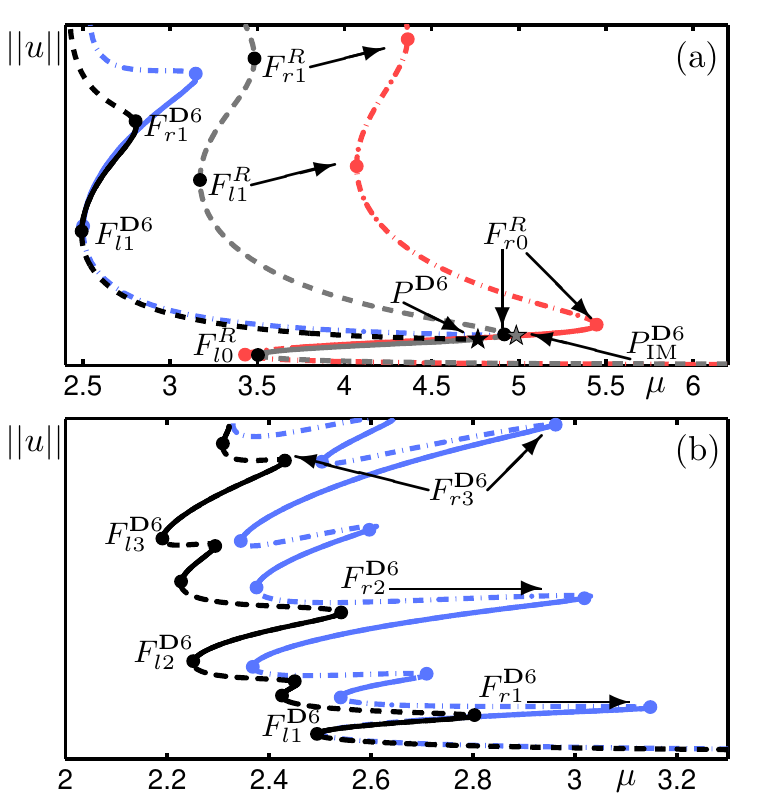}
\caption{Comparison between the solution branches of PDE4 (greyscale)
  and IM (colour). Black curves and grey correspond to the PDE
  approximation and have the line style and labelling conventions as
  \fref{fig:lt4th_snake_1}. Corresponding radial and $\dsix$ solution
  branches from IM are plotted as dash-dot curves in red and blue,
  respectively; stability is not indicated and folds not labelled on
  these branches. (a): Detail close to the spatial-symmetry-breaking
  bifurcation $\pds$; the corresponding bifurcation for IM is labelled
  $\pdsfull$. (b): Detail of the $\dsix$ snaking structure.
  Parameters given at the beginning of \sref{snaking1}.}
\label{fig:lt4th_compare}
\end{figure}
The same radially- and $\dsix$-symmetric branches shown in the
previous section have been computed for IM. Here
we test the accuracy of PDE4 both qualitatively in terms of the types
of solutions produced and their bifurcations, and quantitatively in
terms of the parameter ranges for which the different solution types
persist. We are also interested to see whether the relative ranges of
existence for different types of solution is consistent between
PDE4 and IM.

\Fref{fig:lt4th_compare}(a) and (b) both show detail from
\fref{fig:lt4th_snake_1}(b) with the same curves reproduced for PDE4
with the same line style and labelling conventions. Also plotted (in
colour) are the equivalent curves computed for IM, where the
equivalent of $\pds$ in IM is $\pdsfull$. The first major point to
make is that in terms of the types of solution encountered, the series
of bifurcations encountered and the stability of each branch segment,
there is an exact agreement between PDE4 and IM. Furthermore, the
quantitative agreement on the radial branch is good up until the fold
point $\frr{1}$. Above $\frr{1}$, the radial branch for PDE4 makes a
large excursion away from the IM branch and the branches remain well
separated as the snaking continues; see panel (a). Similarly, for the
$\dsix$ branch the level of agreement is good up to $\fdsr{1}$ above
which PDE4 branch deviates and remains well separated from the IM
branch; see panel (b). The range of existence in $\mu$ for each stable
branch segment on the $\dsix$ branch is greatly under estimated by
PDE4.

We now highlight a key qualitative difference between the bifurcation
diagrams for PDE4 and IM. In  IM, there is a range of
$\mu\in[2.5,3.0]$ for which the stable branch segments corresponding
to $7$-spot, $19$-spot and $37$-spot all overlap. This is not the case
for PDE4, in particular, the branch segments corresponding to
stable $7$-spot and $37$-spot solutions between the fold-pairs
$(\fdsl{1},\fdsr{1})$ and $(\fdsl{3},\fdsr{3})$ do not overlap. This
can be seen by the fact that $\fdsr{3}$ occurs at a smaller
$\mu$-value (indicated by the first thin vertical line in
\fref{fig:lt4th_snake_1}(b)) than $\fdsl{1}$ (indicated by the second
thin vertical line). This organisation of the solutions in parameter
space is qualitatively inconsistent with IM.

\subsubsection{PDE8 results and comparison with IM}\label{sec:8th}

\begin{figure}
\centering
\includegraphics[width=\linewidth]{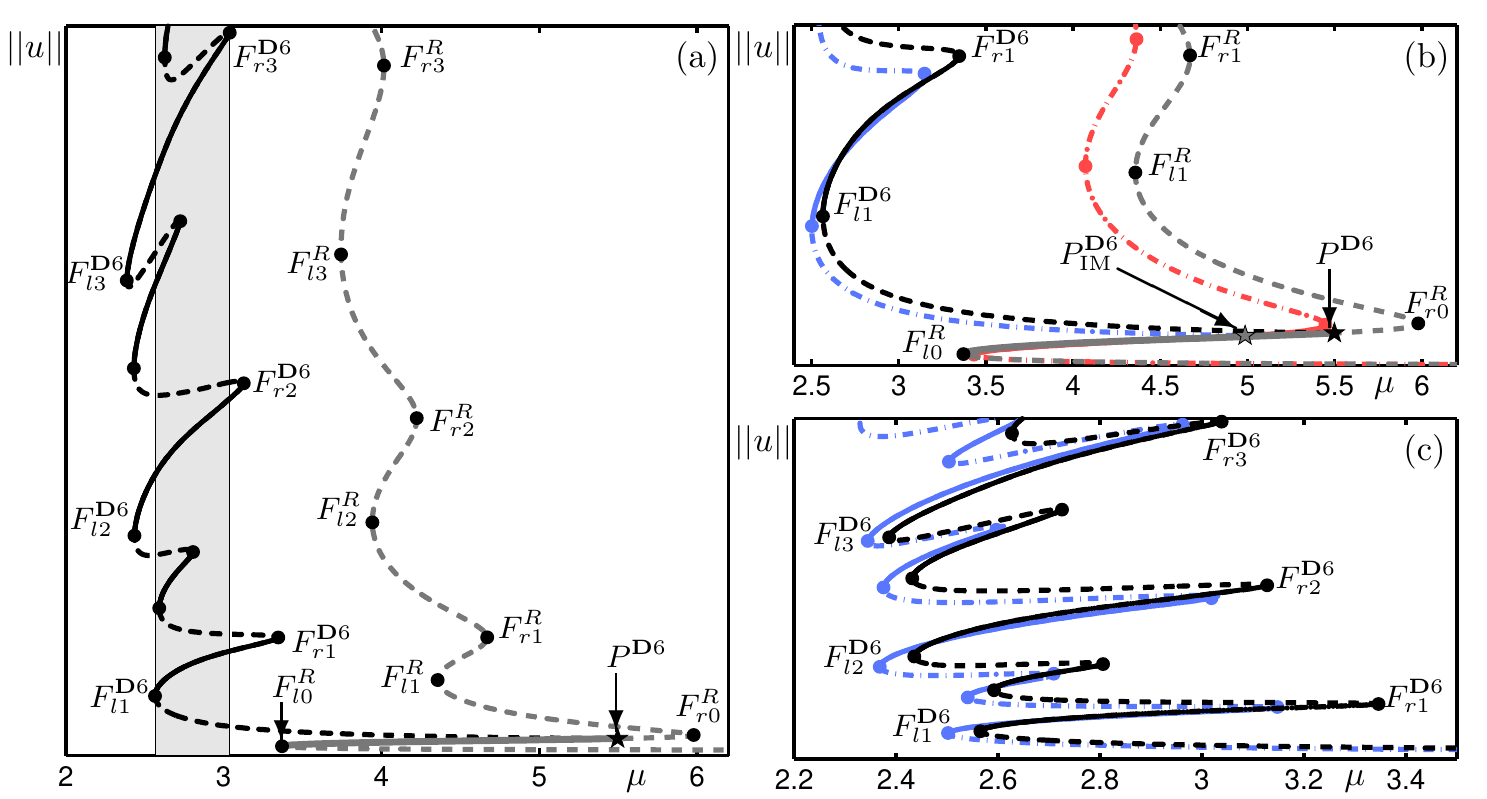}
\caption{Snaking of radially- and $\dsix$-symmetric localised
  solutions in terms of $\mu$ for PDE8 (greyscale) and comparison with
  IM (colour). Black and grey curves correspond to the PDE
  approximation and have the line style and labelling conventions
  as~\fref{fig:lt4th_snake_1}. Corresponding radial and $\dsix$
  solution branches from the full model are plotted as dash-dot curves
  in red and blue, respectively; stability is not indicated and folds
  not labelled on these branches. (a): Global snaking diagram for
  PDE8. (b): Detail close to the spatial-symmetry-breaking bifurcation
  $\pds$; the corresponding bifurcation in the full model is labelled
  $\pdsfull$. (c): Detail of the $\dsix$ snaking structure.
  Parameters given at the beginning of \sref{snaking1}.}
\label{fig:lt8th_snake_1}
\end{figure}

\Fref{fig:lt8th_snake_1}(a) shows branches of both radially- and
$\dsix$-symmetric solutions of PDE8. Globally the bifurcation diagram
is the same as that of PDE4 in terms of the types of solution observed
and the sequence of bifurcation encountered. There is an important
difference between PDE4 and PDE8 in terms of the organisation in
parameter space of the solution branches. For PDE4, there is no
overlap in parameter ranges for which the $7$- and $37$-spot branches
are stable; see the two vertical lines in \fref{fig:lt4th_snake_1} and
note that $\fdsr{3}$ occurs before $\fdsl{1}$ in this case. In the
PDE8 case, as shown in \fref{fig:lt8th_snake_1}(a), there is an
overlap in the parameter ranges as indicated by the grey shaded
region. This organisation of the solution branches in parameter space
is now consistent with the full model as shown in panel (c). Indeed
PDE8 provides better agreement with IM; the branches remain close for
both the radial and the $\dsix$ branches as we move up the snaking
structure as can be seen in panels (b) and (c). We note that when
compared with IM, the upper section of the radial branch occurs at
smaller values of $\mu$ for PDE4 and at larger values of $\mu$ for
PDE8; compare \fref{fig:lt4th_snake_1}(a) with
\fref{fig:lt8th_snake_1}(b).

\subsection{Snaking of $\dtwo$, $\dthree$ and $\dfour$ patterns in IM}\label{sec:numerics_full}\label{snaking2}
In this section we discuss several patterns that possess neither
radial nor $\dsix$ symmetry. For non-radial patterns discussed so far
in this article, see \fref{fig:lt4th_snake_1}(b1)--(b3), the
individual spots lie on a regular hexagonal lattice. However, there
exist an infinite number of stable configurations that conform to the
same lattice spacing but without the full $\dsix$ symmetry; here we
present two such examples. We also present one further example of a
stable configuration that does not conform to the regular hexagonal
lattice.

\subsubsection{Three-spot pattern ($\dthree$)}\label{sec:trirad}
\begin{figure}
\centering
\includegraphics[width=0.8\textwidth]{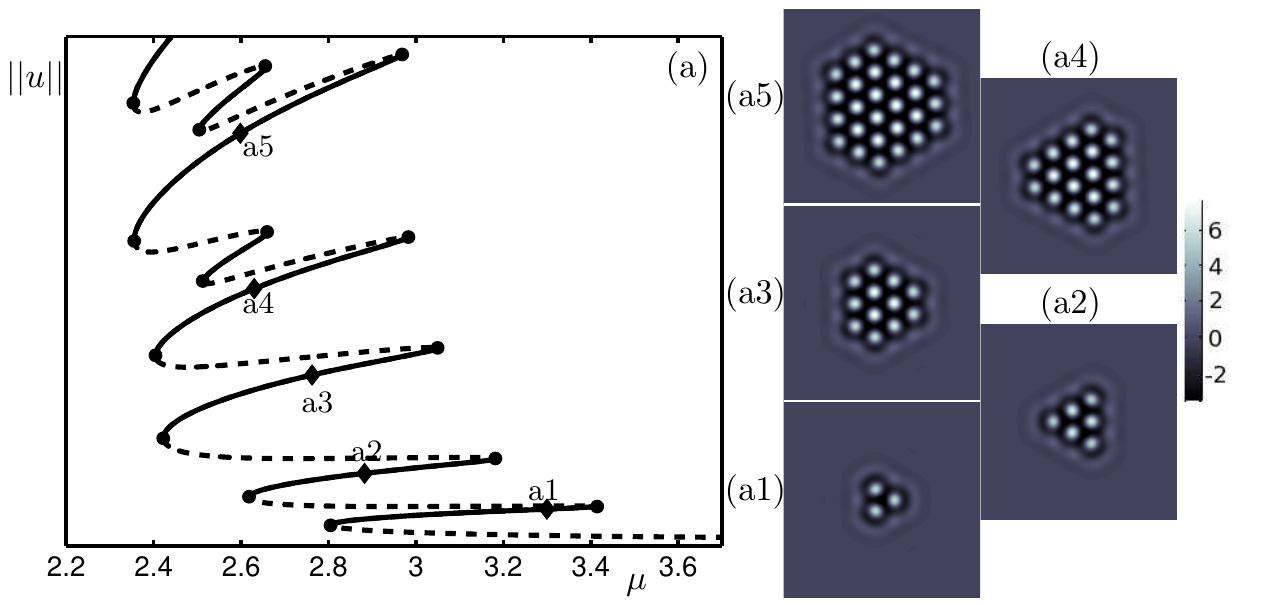}
\caption{Snaking of patterns with $\dthree$ symmetry built around a
  triangular-organised, three-spot solution. Line style conventions
  for stability as in~\fref{fig:lt4th_snake_1}.  Parameters for IM as
  given at the beginning of \sref{snaking1}.}
\label{fig:full_trirad_1}
\end{figure}

\Fref{fig:full_trirad_1} shows the snaking of $\dthree$-symmetric
patterns about a three-spot solution. The unstable branch that enters
panel (a) in the bottom-right-hand corner reconnects to the trivial
state $u=0$ at the Turing instability (not shown). The panels
(a1)--(a5) show stable solutions on the first five full excursions in
$\mu$ of the snaking structure. The existence of three-spot (see panel
(a1)) and twelve-spot (see panel (a3)) solutions for PDE4 with
connectivity given by \eqreff{eq:what_2d_4th} was shown
in~\cite{laing-troy:03}. Here we have shown that these solutions exist
in IM and that they form part of a larger snaking structure.

\subsubsection{Two-spot pattern ($\dtwo$)}\label{sec:twospot}
\begin{figure}
\centering
\includegraphics[width=0.8\textwidth]{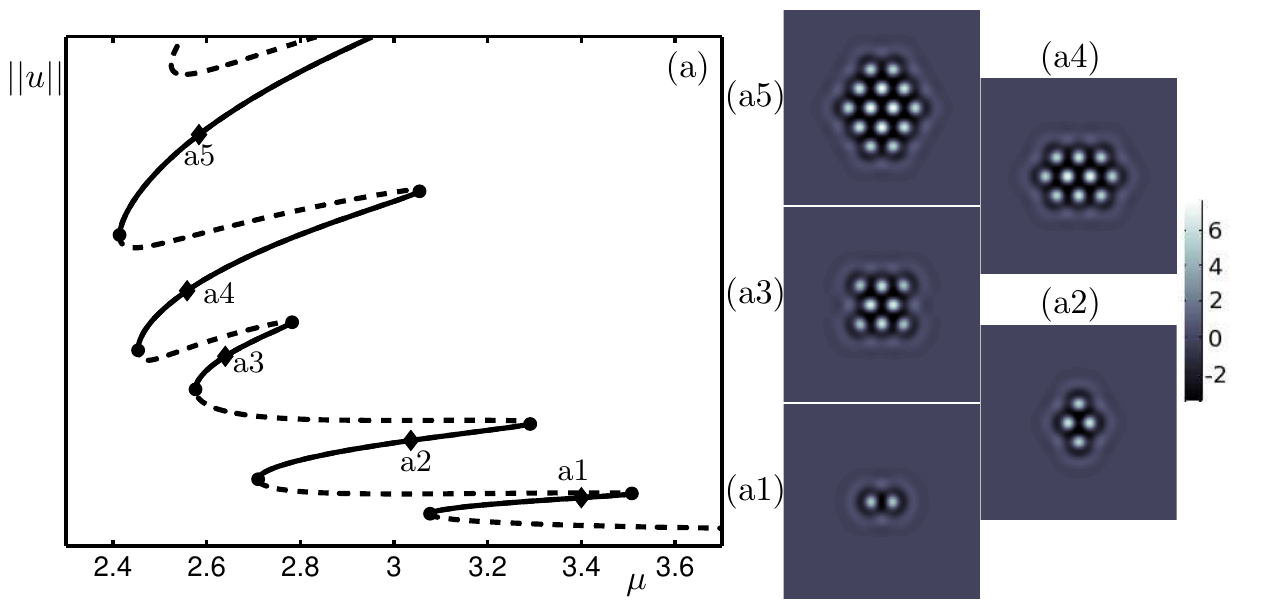}
\caption{Snaking of patterns with $\dtwo$ symmetry built around a
  two-spot solution. Line style conventions for stability as
  in~\fref{fig:lt4th_snake_1}.  Parameters for IM given at the
  beginning of \sref{snaking1}.}
\label{fig:full_twospot_1}
\end{figure}

Similarly, there is an unstable two-spot solution that connects to the
trivial state $u=0$ at the Turing instability (not
shown). \Fref{fig:full_twospot_1}(a) shows that this solution also
undergoes a sequence of fold bifurcations giving rise to larger
$\dtwo$-symmetric patterns. We note that the spacing between the spots
in these patterns still conforms to the regular hexagonal lattice. The
panels (a1)--(a5) show solutions on the first five stable branch
segments moving up the snaking structure; we note that the pattern
(a3) is on an intermediate branch that does not make a full excursion
in $\mu$.

\subsubsection{Quincunx pattern ($\dfour$)}\label{sec:quinc}
\begin{figure}
\centering
\includegraphics[width=0.65\textwidth]{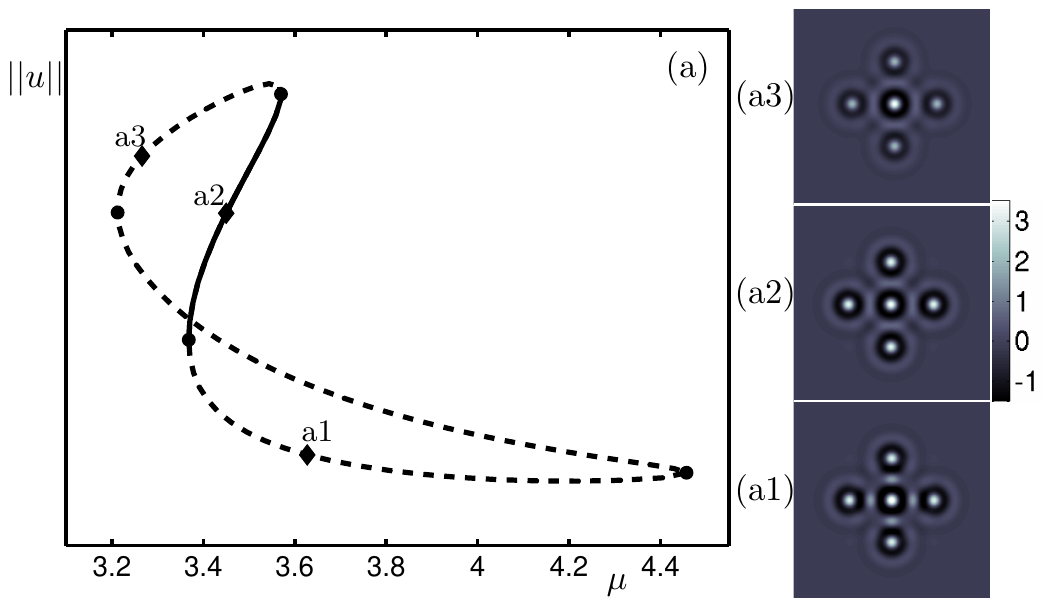}
\caption{Isola of five-spot solution organised in a quincunx
  configuration. Line style conventions for stability as
  in~\fref{fig:lt4th_snake_1}. Parameters for IM given at the
  beginning of \sref{snaking1}.}
\label{fig:full_quinc_1}
\end{figure}

We show in \fref{fig:full_quinc_1}(a2) a stable configuration that
lies on a square lattice but with a spacing between the spot peaks
that is double that of the hexagonal-lattice solutions encountered
thus far. The solution is formed by five spots that interact at the
first excitatory peak away from $0$ in the connectivity function
(\ref{eq:intModelW}) as shown in \fref{fig:full_conn_sig}(d). The
configuration of four spots forming a square with an additional spot
in the center is typically referred to as a quincunx pattern that is
found, for example, on dice and dominoes. As shown in
\fref{fig:full_twospot_1}(a) these solutions exist on an isola in
parameter space where other branches, see panels (a1) and (a3), are
unstable. Although the pattern does not undergo snaking-type behaviour
it may be possible to construct larger patterns on the double-spaced
square lattice.

\section{Discussion}
\subsection{Summary}
This paper explores patterns of localised activity in the neural field
equation posed on the Euclidean plane with a smooth firing rate
function.  The choice of connectivity function is an important factor
in determining whether, the localised behaviour found is restricted to
individual spots, or whether multiple interacting spots can form
coherent localised patterns. In~\cite{faye-rankin-etal:12b} localised
states were studied in a model with a radially-symmetric wizard-hat
connectivity function describing local excitation and lateral
inhibitions in the Euclidean plane. When spot solutions were tracked
using numerical continuation no snaking behaviour was observed,
i.e. the only steady states found consisted of a single spot. The
radially symmetric connectivity function studied here and shown in
\fref{fig:full_conn_sig}(d)) features local excitation, lateral
inhibition and long-range bands of excitation that decay with
distance. Indeed, the distance between excitation peaks fixes a
spatial scale that allows for regular spatial interactions and the
formation of larger patterns of activity. In~\cite{laing-troy:03} it
was shown that multiple-spot patterns could be obtained all with the
common property of the peaks lying on a regular hexagonal lattice with
spacing determined by the connectivity function. One of the main aims
in the present article was to show how these solutions are connected
in parameter space and how patterns with varying spatial extent grow
via the mechanism of homoclinic snaking.

The results presented in~\cite{laing-troy:03} relied on working with
an approximated connectivity function (see~\fref{fig:ltpde_conn}(a)
and (b)) that allowed for solutions the full integral neural field
equation to be studied in an equivalent fourth-order PDE . The initial
parts of the results section are concerned with the relative agreement
between solutions to the integral model and equivalent PDE
formulations with approximated connectivity functions. We pursue the
problem numerically by investigating the level of agreement in terms
of an entire bifurcation diagram rather than comparing individual
solutions at fixed parameter values. We compared both radially
symmetric and $\dsix$-symmetric solution branches and found that the
qualitative difference of the zero-mode in the Fourier domain for the
approximated connectivity used in the fourth-order model
(see~\fref{fig:ltpde_conn}(b)) led to significant discrepancies in the
existence ranges of solution branches, in particular for solutions
with a larger spatial extent. We demonstrated that the improved
approximation of the connectivity function shown
in~\fref{fig:ltpde_conn}(c) and (d), leading to an eighth-order PDE,
provides a better agreement across the full bifurcation diagram. In
particular, the eigth-order model captures the key feature of there
being a specific parameter range in which multiple solutions coexist,
each solution with a different spatial extent. It was not possible to
capture this feature with the fourth-order approximation. We conclude
that, although equivalent PDE formulations have proved to be a useful
tool for the study of neural fields, it is important to ensure close
agreement between the connectivity functions in the Fourier
domain. Increasing the order of the PDEs used allows for improvements
in this agreement.  We believe that, while converting the integral
formulation to higher-order PDEs could be useful in analytical
studies, numerical calculations of these systems should be approached
without resorting to PDE formulation where possible. In passing, we
point out that the methodology proposed here for the integral model is
applicable to inhomogeneous synaptic kernels, provided that the
convolution structure of the integral is preserved. Furthermore, we
point out that here we have used the standard Newton-GMRES method
mainly for its simplicity, but more sophisticated choices are also
possible~\cite{kelley:03}.

Having investigated radial and $\dsix$ solutions with PDE formulations
and compared the results with the full integral model, we proceeded to
give an account of other types of solutions that, when path-followed
with numerical continuation lead to patterns with different underlying
symmetry properties. We worked with the full integral model and showed
that patterns with $\dtwo$ and $\dthree$ symmetry also give rise to
snaking behaviour, generating spatial patterns with variable spatial
extent. All the solutions of this type, including those shown earlier
with $\dsix$-symmetry, have the common feature of the individual spots
lying on a regular hexagonal lattice. Furthermore, these solutions of
variable spatial extent exist within roughly the same ranges of the
parameter $\mu$. As the snaking diagram is ascended, new spots are
glued to long edges of the pattern in a regular fashion. We note the
existence of an arbitrary number of intermediate branches not shown in
our bifurcation diagrams where symmetries can be broken via the
(simultaneous) addition or subtraction of one (or more spots). We
expect these solutions to lie on intermediate branches that are stable
over ranges smaller than the full excursions in the main snaking
structure.

\begin{figure}
\centering
\includegraphics[width=\linewidth]{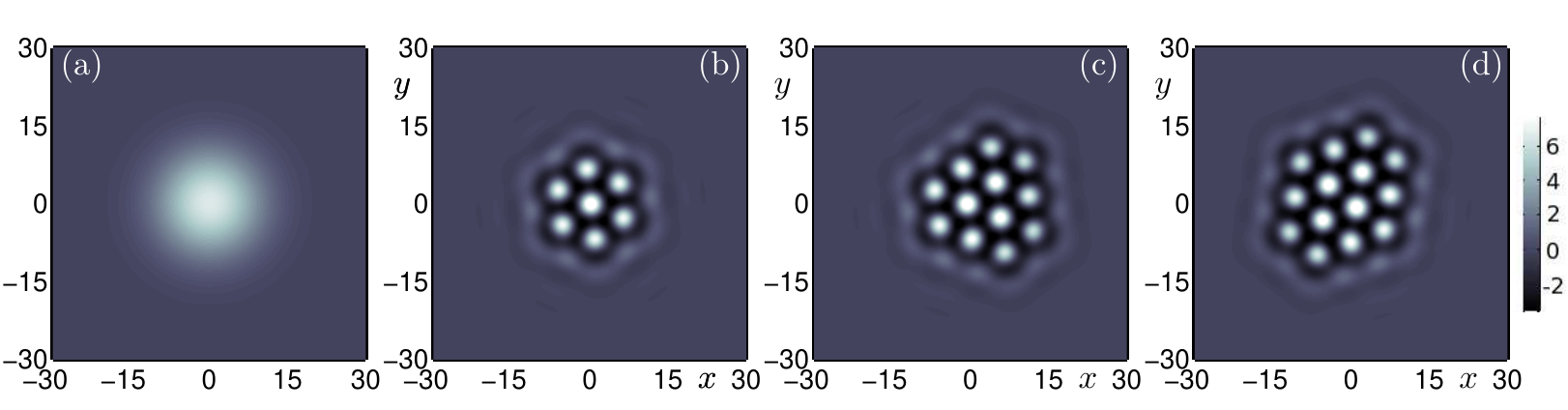}
\caption{Model simulations with input given by \eqreff{eq:extInput}
  with $\alpha=\beta=1$. The simulations are performed with
  $\mu=2.4$. (a) The Gaussian-bump input plotted with $\sigma=10$ and
  $G_0=6$ ($G_0=1.5$ used in the simulations). (b): Case with
  $\sigma=9.0$, the model converges to a $7$s-spot solution that is a
  modified version of the state shown in~\fref{fig:lt4th_snake_1}(b1)
  (c): Case with $\sigma=9.5$, the model converges to a $12$s-spot
  solutions that is a modified version of the state shown
  in~\fref{fig:full_trirad_1}(a3) (d): Case with $\sigma=10.0$, the
  model converges to a $14$s-spot solutions that is a modified version
  of the state shown in~\fref{fig:full_twospot_1}(a5).  Other model
  parameters (for IM) given at the beginning of \sref{snaking1}. }
\label{fig:input_th}
\end{figure}

\subsection{Localised patterns with input}\label{sec:discuss_inp}
The numerical bifurcation analysis presented in this paper opens up
the possibility to investigate the spread of cortical activity with
inputs.  An important principle in studying the neural field equation
is that weak inputs to the equations should drive the system to states
that are already solutions to the underlying equations without
input. The bifurcation study presented here allows us to identify the
types of solution that we may expect to encounter for localised inputs
and the relevant parameter regime in which they occur. Two criteria
need to be satisfied when identifying a suitable operating regime, 1)
the model should only produce the trivial homogeneous state $u=0$
before an input is introduced and 2) when an input is introduced, the
model should be driven to one of the underlying non-trivial solutions.
We have shown that, for the full integral model, there is an
accumulation of fold bifurcations at around $\mu=2.4$ representing the
first point for which localised patterns can be observed. Both
criteria are satisfied when the model is operated just before these
fold points. Introducing a weak input the system can be driven to
states that have a spatial extent corresponding to that of the
input. \Fref{fig:input_th}(a) shows the profile of the input used in a
series of simulations that were initiated with the $u=0$ state plus a
small random perturbation across the entire spatial
domain. \Frefs{fig:input_th}(b), (c) and (d) show the result of three
$150$ time unit simulations, each with different spatial extent for
the input. In each case, the trivial state $u=0$ is no longer stable
and the system naturally selects one of the solutions described in the
early bifurcation analysis. As the spatial extent of the input
increases, the size of the pattern selected by the model increases and
this is an important consequence of the corresponding solutions
existing as part of a snaking structure. The computational framework
presented in this article will allow for the relationship between
model inputs and spatially localised patterns to be investigated in
future work.

\subsection{Patterns on other lattices and parametric forcing}
\begin{figure}
\centering
\includegraphics[width=\linewidth]{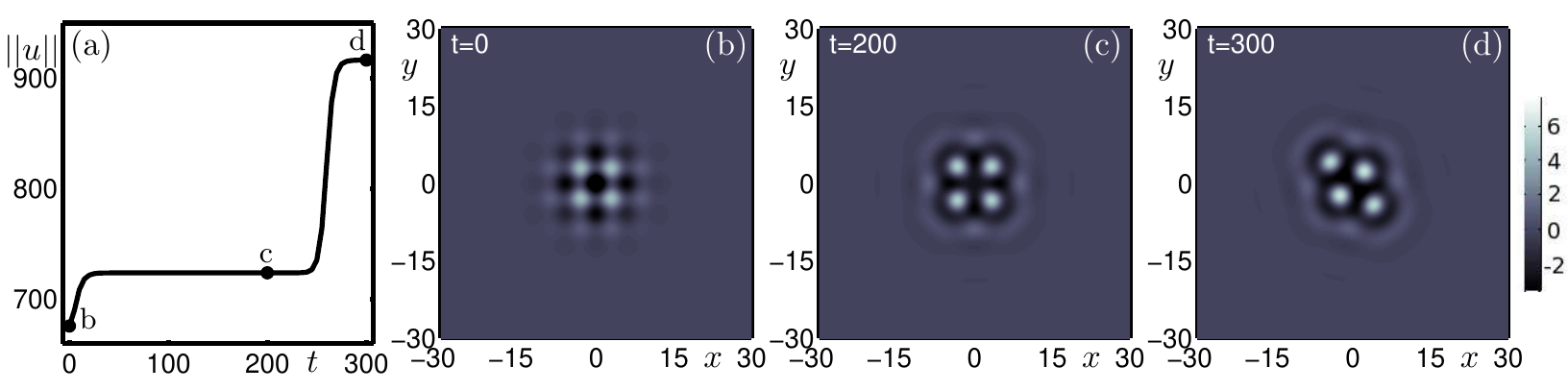}
\caption{Symmetry breaking of weakly unstable solution on a square
  lattice at $\mu=3.2$. (a): Evolution of L2 norm. (b): Initial
  condition given by \eqreff{eq:d4ic}. (c): Weakly unstable
  solution. (d): Stable solution. Parameters (for IM) given at the
  beginning of \sref{snaking1}.}
\label{fig:full_d4_th}
\end{figure}

The multi-spot solutions described in this article all have the common
feature of the activated peaks falling onto a regular hexagonal
lattice. Indeed, this has been found to be the default way for the
radial symmetry to be broken in pattern forming systems, notably the
archetypal Swift-Hohenberg equation~\cite{lloyd-sandstede-etal:08}. We
also investigated whether it is possible to find other states that do
not conform to the regular hexagonal lattice. In
\fref{fig:full_quinc_1} an example of a solution with $\dfour$
symmetry was shown that consists of five spots interacting at double
the standard separation between excitation peaks. We found this
solution to exist on an isola and not undergo snaking so as to find
larger patterns tiling the plane with the same spacing. We also
attempted to converge solutions with $\dfour$ symmetry that have the
regular spacing between peaks. A suitable initial condition to find
such solutions is given by~\eqreff{eq:d4ic} chosen such that there is
a depression at $(x,y)=(0,0)$ and the surrounding square-lattice
pattern decays away from the origin, see \fref{fig:full_d4_th}(b).  In
an appropriate parameter range these states appear to converge to
stable patterns. However, we found the patterns to be weakly unstable,
finally converging to a pattern on a hexagonal lattice after a long
transient. \Fref{fig:full_d4_th}(a) shows a time-course of the L2 norm
from a simulation with the initial condition shown in panel (b). The
model reaches the apparently stable $\dfour$ configuration after
approximately $30$ time units (panel (c)) before finally converging
after $230$ time units to the $\dtwo$-symmetric pattern (panel (d))
that was previously identified in \fref{fig:full_twospot_1}(a2). It
would be possible to stabilise such weakly unstable solutions with the
use of parametric forcing, by introducing small modulations that
encourage interactions on a fixed lattice in cortical space. In the
neural field equations this is typically referred to as inhomogeneous
neural media; travelling waves, travelling fronts, periodic patterns
and pulsating fronts have been studied with such
modulations~\cite{bressloff:01,coombes-venkov-etal:07,coombes-laing:11}. The
study of localised states in this context would be an interesting
future direction. Furthermore, stable quasi-periodic patterns have
been obtained in the Swift-Hohenberg equation through parametric
forcing~\cite{iooss-rucklidge:10}.  The question of introducing
orientation-preference tilings on square and hexagonal lattices was
addressed in~\cite{baker-cowan:09}. However, one is restricted in the
number of orientations that can be equally represented with such
tilings, (two and three, respectively). Parametric forcing on a
quasi-periodic lattice is of particular interest in the neural field
equations as this could allow for near-continuous representations of
features in a model without an abstracted feature space.

\section{Conclusions}
The organisation in parameter space of localised structures consisting
of multiple spots has been revealed for the first time in planar
neural field equations. In order to find such behaviour, one must
choose a connectivity function with an excitatory peak away from the
origin that fixes a regular spatial scale of interactions between
spots. As localised solutions are path-followed using numerical
continuation we find that these structures grow in a series of fold
bifurcations through the mechanism of homoclinic snaking that has been
well-studied in the Swift-Hohenberg equation. A numerical strategy has
been proposed to perform a numerical bifurcation analysis without
resorting to a PDE formulation, but taking advantage of matrix-free
Newton-Krylov nonlinear solvers combined with a pseudospectral
evaluation of the right-hand side. The novel application of these
methods to the neural field equations allowed for numerical
continuation to be applied to the full integral form of the
model. Previous studies in 2D have relied exclusively on PDE
approximations of the connectivity functions; here we demonstrated
that these approximations can give a very close agreement with the
full integral model if a sufficiently high-order approximation is
taken. The numerical schemes presented here will allow for future
studies of the neural field equations to use connectivity functions
defined either directly in the real domain or the Fourier domain
without recourse to PDE methods, provided that the sigmoidal firing
rate be smooth and that the integral formulation can be expressed as a
convolution (this extends also to inhomogeneous firing rates).

The neural field studied in the present paper can be considered as a
model of the visual cortex and the localised patterns studied without
inputs can be related to visual hallucinations that can be localised in
the visual field~\cite{siegel:77}. Furthermore, we have shown that the
localised states computed in our bifurcation analysis are exactly the
types of solutions selected by the model in the presence of weak
inputs. The persistence of these localised structures in the presence
of a model input is new.  This future direction will be of particular
interest for the study of localised patterns of activity that have
been observed in the primary visual
cortex~\cite{chavane-sharon-etal:11} with localised visual input.

\section*{Acknowledgments}
The research leading to these results has received funding from ERC
grant 227747 (NERVI) {\bf(JR,JB)}, Marie-Curie Facets-ITN grant 237955
{\bf(JB)}, and the BrainScaleS project grant 269921 {\bf(JB)}. {\bf
  DL} gratefully acknowledges partial support of the UK Engineering
and Physical Sciences Research Council grant EP/H05040X/1 (Nucleation
of Ferrosolitons and Localised Ferro-patterns).

\appendix
\setlength{\belowdisplayskip}{0pt} \setlength{\belowdisplayshortskip}{0pt}
\setlength{\abovedisplayskip}{0pt} \setlength{\abovedisplayshortskip}{0pt}
\section{Analytic expressions for initial conditions}
\Fref{fig:full_th}(a) shows the initial
condition given by
\begin{equation}\label{eq:radic}
u(x,y)=A\exp \bigg(-\frac{x^2+y^2}{L}\bigg),
\end{equation}
with $A=6$ and $L=5.77$; subsequent panels (b) and (c) show a
transient state after $1$ time unit and the stable steady state after
$15$ time units. \Fref{fig:full_th}(d) shows the initial condition
given by
\begin{equation}\label{eq:hexic}
u(x,y)=A\exp \bigg(-\frac{x^2+y^2}{L}\bigg)\left[\cos(x)+\cos\left(\frac{1}{2}x+\frac{\sqrt{3}}{2}y\right) +\cos\left(-\frac{1}{2}x+\frac{\sqrt{3}}{2}y\right) \right],
\end{equation}
with $A=2$ and $L=100$; subsequent panels (e) and (f) show a
transient state after $1$ time unit and the stable steady state after
$15$ time units. \Fref{fig:full_d4_th}(b) shows the initial condition given by
\begin{equation}\label{eq:d4ic}
u(x,y)=2\exp \bigg(-\frac{x^2+y^2}{L}\bigg)(-\cos x-\sin y),
\end{equation}
with $A=2$ and $L=65$; subsequent panels (c) and (d) show transient
states after $200$ and $300$ time units, respectively. The $L_2$-norm
is plotted over this time course in panel (a).


\begin{thebibliography}{10}

\bibitem{amari:75}
{\sc S.~Amari}, {\em Homogeneous nets of neuron-like elements}, Biological
  Cybernetics, 17 (1975), pp.~211--220.

\bibitem{amari:77}
{\sc S.-I. Amari}, {\em Dynamics of pattern formation in lateral-inhibition
  type neural fields}, Biological Cybernetics, 27 (1977), pp.~77--87.

\bibitem{avitabile-lloyd-etal:10}
{\sc D.~Avitabile, D.~Lloyd, J.~Burke, E.~Knobloch, and B.~Sandstede}, {\em To
  snake or not to snake in the planar swift-hohenberg equation}, SIAM J. Appl.
  Dyn. Syst., 9 (2010), pp.~704--733.

\bibitem{baker-cowan:09}
{\sc T.~Baker and J.~Cowan}, {\em Spontaneous pattern formation and pinning in
  the primary visual cortex}, Journal of Physiology-Paris, 103 (2009),
  pp.~52--68.

\bibitem{beck-knobloch-etal:09}
{\sc M.~Beck, J.~Knobloch, D.~Lloyd, B.~Sandstede, and T.~Wagenknecht}, {\em
  Snakes, ladders, and isolas of localised patterns}, SIAM J. Math. Anal, 41
  (2009), pp.~936--972.

\bibitem{ben-yishai-bar-or-etal:95}
{\sc R.~Ben-Yishai, R.~Bar-Or, and H.~Sompolinsky}, {\em Theory of orientation
  tuning in visual cortex}, Proceedings of the National Academy of Sciences, 92
  (1995), pp.~3844--3848.

\bibitem{bressloff:01}
{\sc P.~Bressloff}, {\em Traveling fronts and wave propagation failure in an
  inhomogeneous neural network}, Physica D: Nonlinear Phenomena, 155 (2001),
  pp.~83--100.

\bibitem{bressloff:12}
\leavevmode\vrule height 2pt depth -1.6pt width 23pt, {\em Spatiotemporal
  dynamics of continuum neural fields}, Journal of Physics A: Mathematical and
  Theoretical, 45 (2012).

\bibitem{bressloff-cowan-etal:01}
{\sc P.~Bressloff, J.~Cowan, M.~Golubitsky, P.~Thomas, and M.~Wiener}, {\em
  Geometric visual hallucinations, euclidean symmetry and the functional
  architecture of striate cortex}, Phil. Trans. R. Soc. Lond. B, 306 (2001),
  pp.~299--330.

\bibitem{bressloff-kilpatrick:11}
{\sc P.~Bressloff and Z.~Kilpatrick}, {\em Two-dimensional bumps in piecewise
  smooth neural fields with synaptic depression}, SIAM Journal on Applied
  Mathematics, 71 (2011), pp.~379--408.

\bibitem{burke-knobloch:07}
{\sc J.~Burke and E.~Knobloch}, {\em Homoclinic snaking: structure and
  stability}, Chaos, 17 (2007), p.~7102.

\bibitem{burke-knobloch:07c}
\leavevmode\vrule height 2pt depth -1.6pt width 23pt, {\em Snakes and ladders:
  localized states in the swift--hohenberg equation}, Physics Letters A, 360
  (2007), pp.~681--688.

\bibitem{buzas-eysel-etal:01}
{\sc P.~Buz{\'a}s, U.~Eysel, P.~Adorj{\'a}n, and Z.~Kisv{\'a}rday}, {\em Axonal
  topography of cortical basket cells in relation to orientation, direction,
  and ocular dominance maps}, The Journal of comparative neurology, 437 (2001),
  pp.~259--285.

\bibitem{canuto-hussaini:06}
{\sc C.~Canuto, M.~Hussaini, A.~Quarteroni, and T.~Zang}, {\em Spectral
  Methods: Fundamentals in Single Domains}, Springer-Verlag, Berlin, 2006.

\bibitem{champneys-sandstede07}
{\sc A.~Champneys and B.~Sandstede}, {\em Numerical computation of coherent
  structures}, in Numerical Continuation Methods for Dynamical Systems,
  Springer, 2007, pp.~331--358.

\bibitem{chapman-kozyreff:09}
{\sc S.~Chapman and G.~Kozyreff}, {\em Exponential asymptotics of localized
  patterns and snaking bifurcation diagrams.}, Physica~D, 238 (2009),
  pp.~319--354.

\bibitem{chavane-sharon-etal:11}
{\sc F.~Chavane, D.~Sharon, D.~Jancke, O.~Marre, Y.~Fr{\'e}gnac, and
  A.~Grinvald}, {\em Lateral spread of orientation selectivity in v1 is
  controlled by intracortical cooperativity}, Frontiers in Systems
  Neuroscience, 5 (2011).

\bibitem{colby-duhamel-etal:95}
{\sc C.~Colby, J.~Duhamel, and M.~Goldberg}, {\em Oculocentric spatial
  representation in parietal cortex}, Cereb. Cortex, 5 (1995), pp.~470--481.

\bibitem{coombes:05}
{\sc S.~Coombes}, {\em Waves, bumps, and patterns in neural fields theories},
  Biological Cybernetics, 93 (2005), pp.~91--108.

\bibitem{coombes-laing:11}
{\sc S.~Coombes and C.~Laing}, {\em Pulsating fronts in periodically modulated
  neural field models}, Physical Review E, 83 (2011), p.~011912.

\bibitem{coombes-lord-etal:03}
{\sc S.~Coombes, G.~Lord, and M.~Owen}, {\em Waves and bumps in neuronal
  networks with axo-dendritic synaptic interactions}, Physica D: Nonlinear
  Phenomena, 178 (2003), pp.~219--241.

\bibitem{coombes-schmidt-etal:13}
{\sc S.~Coombes, H.~Schmidt, and D.~Avitabile}, {\em Neural Field Theory},
  Springer, 2013, ch.~Spots: Breathing, drifting and scattering in a neural
  field model.

\bibitem{coombes-schmidt-etal:12}
{\sc S.~Coombes, H.~Schmidt, and I.~Bojak}, {\em Interface dynamics in planar
  neural field models}, Journal of mathematical neuroscience,  (2012).

\bibitem{coombes-venkov-etal:07}
{\sc S.~Coombes, N.~Venkov, L.~Shiau, I.~Bojak, D.~Liley, and C.~Laing}, {\em
  Modeling electrocortical activity through improved local approximations of
  integral neural field equations}, Physical Review E, 76 (2007), p.~051901.

\bibitem{coullet-riera-etal:00}
{\sc P.~Coullet, C.~Riera, and C.~Tresser}, {\em Stable static localized
  structures in one dimension}, Phys. Rev. Lett., 84 (2000), pp.~3069--3072.

\bibitem{elvin-laing-etal:10}
{\sc A.~Elvin, C.~Laing, R.~McLachlan, and M.~Roberts}, {\em Exploiting the
  hamiltonian structure of a neural field model}, Physica D: Nonlinear
  Phenomena, 239 (2010), pp.~537--546.

\bibitem{ermentrout:98}
{\sc G.~Ermentrout}, {\em Neural networks as spatio-temporal pattern-forming
  systems}, Reports on Progress in Physics, 61 (1998), pp.~353--430.

\bibitem{ermentrout-cowan:79b}
{\sc G.~Ermentrout and J.~Cowan}, {\em A mathematical theory of visual
  hallucination patterns}, Biological Cybernetics, 34 (1979), pp.~137--150.

\bibitem{faye-rankin-etal:12}
{\sc G.~Faye, J.~Rankin, and P.~Chossat}, {\em Localized states in an unbounded
  neural field equation with smooth firing rate function: a multi-parameter
  analysis}, Journal of Mathematical Biology,  (2012).

\bibitem{faye-rankin-etal:12b}
{\sc G.~Faye, J.~Rankin, and D.~Lloyd}, {\em Localized radial bumps of a neural
  field equation on the euclidean plane and the poincar{\'e} disk},
  Nonlinearity, (accepted) (2012).

\bibitem{folias-bressloff:04}
{\sc S.~Folias and P.~Bressloff}, {\em Breathing pulses in an excitatory neural
  network}, SIAM Journal on Applied Dynamical Systems, 3 (2004), pp.~378--407.

\bibitem{folias-bressloff:05}
\leavevmode\vrule height 2pt depth -1.6pt width 23pt, {\em Breathers in
  two-dimensional neural media}, Physical Review Letters, 95 (2005), p.~208107.

\bibitem{funahashi-bruce-etal:89}
{\sc S.~Funahashi, C.~Bruce, and P.~Goldman-Rakic}, {\em Mnemonic coding of
  visual space in the monkey's dorsolateral prefrontal cortex}, J.
  Neurophysiol., 61 (1989), pp.~331--349.

\bibitem{gutkin-ermentrout-etal:00}
{\sc B.~Gutkin, G.~B.~Ermentrout, and J.~O'Sullivan}, {\em Layer 3 patchy
  recurrent excitatory connections may determine the spatial organization of
  sustained activity in the primate prefrontal cortex}, Neurocomputing, 32
  (2000), pp.~391--400.

\bibitem{hansel-sompolinsky:97}
{\sc D.~Hansel and H.~Sompolinsky}, {\em Modeling feature selectivity in local
  cortical circuits}, Methods of neuronal modeling,  (1997), pp.~499--567.

\bibitem{haudin-rojas-etal:11}
{\sc F.~Haudin, R.~Rojas, U.~Bortolozzo, S.~Residori, and M.~Clerc}, {\em
  Homoclinic snaking of localized patterns in a spatially forced system},
  Physical Review Letters, 107 (2011), p.~264101.

\bibitem{iooss-rucklidge:10}
{\sc G.~Iooss and A.~Rucklidge}, {\em On the existence of quasipattern
  solutions of the swift--hohenberg equation}, Journal of Nonlinear Science, 20
  (2010), pp.~361--394.

\bibitem{kelley:95}
{\sc C.~Kelley}, {\em Iterative methods for linear and nonlinear equations},
  Society for Industrial and Applied Mathematics, 1995.

\bibitem{kelley:03}
\leavevmode\vrule height 2pt depth -1.6pt width 23pt, {\em Solving nonlinear
  equations with Newton's method}, Society for Industrial and Applied
  Mathematics, 2003.

\bibitem{knobloch:08}
{\sc E.~Knobloch}, {\em Spatially localized structures in dissipative systems:
  open problems}, Nonlinearity, 21 (2008), pp.~T45--T60.

\bibitem{laing:05}
{\sc C.~Laing}, {\em Spiral waves in nonlocal equations}, SIAM Journal on
  Applied Dynamical Systems, 4 (2005), pp.~588--606.

\bibitem{laing:13}
\leavevmode\vrule height 2pt depth -1.6pt width 23pt, {\em Neural Field
  Theory}, Springer, 2013, ch.~PDE Methods for Two-Dimensional Neural Fields.

\bibitem{laing-troy:03}
{\sc C.~Laing and W.~Troy}, {\em {PDE} methods for nonlocal models}, SIAM
  Journal on Applied Dynamical Systems, 2 (2003), pp.~487--516.

\bibitem{laing-troy-etal:02}
{\sc C.~Laing, W.~Troy, B.~Gutkin, and G.~Ermentrout}, {\em Multiple bumps in a
  neuronal model of working memory}, SIAM J. Appl. Math., 63 (2002),
  pp.~62--97.

\bibitem{lloyd-sandstede:09}
{\sc D.~Lloyd and B.~Sandstede}, {\em Localized radial solutions of the
  {S}wift-{H}ohenberg equation}, Nonlinearity, 22 (2009), pp.~485--524.

\bibitem{lloyd-sandstede-etal:08}
{\sc D.~Lloyd, B.~Sandstede, D.~Avitabile, and A.~Champneys}, {\em Localized
  hexagon patterns of the planar {S}wift--{H}ohenberg equation}, SIAM J. Appl.
  Dynam. Syst., 7 (2008), pp.~1049--1100.

\bibitem{mccalla-sandstede:10}
{\sc S.~McCalla and B.~Sandstede}, {\em Snaking of radial solutions of the
  multi-dimensional swift-hohenberg equation: A numerical study}, Physica D:
  Nonlinear Phenomena, 239 (2010), pp.~1581--1592.

\bibitem{nickolls-buck-etal:08}
{\sc I.~Nickolls, J.and~Buck, M.~Garland, and K.~Skadron}, {\em Scalable
  parallel programming with cuda}, Queue, 6 (2008), pp.~40--53.

\bibitem{niebur-wortgotter:94}
{\sc E.~Niebur and F.~W{\"o}rg{\"o}tter}, {\em Design principles of columnar
  organization in visual cortex}, Neural Computation, 6 (1994), pp.~602--614.

\bibitem{owen-laing-etal:07}
{\sc M.~Owen, C.~Laing, and S.~Coombes}, {\em Bumps and rings in a
  two-dimensional neural field: splitting and rotational instabilities}, {New
  Journal of Physics}, 9 (2007), pp.~378--401.

\bibitem{sakaguchi-brand:96}
{\sc H.~Sakaguchi and H.~Brand}, {\em Stable localized solutions of arbitrary
  length for the quintic {S}wift-{H}ohenberg equation}, Physica~D, 97 (1996),
  pp.~274--285.

\bibitem{sakaguchi-brand:98}
\leavevmode\vrule height 2pt depth -1.6pt width 23pt, {\em Localized patterns
  for the quintic complex {S}wift-{H}ohenberg equation}, Physica~D, 117 (1998),
  pp.~95--105.

\bibitem{schneider-gibsol-etal:10}
{\sc T.~Schneider, J.~Gibson, and J.~Burke}, {\em Snakes and ladders: localized
  solutions of plane couette flow}, Physical review letters, 104 (2010),
  p.~104501.

\bibitem{siegel:77}
{\sc R.~Siegel}, {\em Hallucinations.}, Scientific American,  (1977).

\bibitem{taylor:99}
{\sc J.~Taylor}, {\em Towards the networks of the brain: from brain imaging to
  consciousness}, Neural Networks, 12 (1999), pp.~943--959.

\bibitem{trefethen:00}
{\sc L.~Trefethen}, {\em Spectral methods in MATLAB}, vol.~10, Society for
  Industrial and Applied Mathematics, 2000.

\bibitem{werner-richter:01}
{\sc H.~Werner and T.~Richter}, {\em Circular stationary solutions in
  two-dimensional neural fields}, Biological Cybernetics, 85 (2001),
  pp.~211--217.

\bibitem{wilson-cowan:72}
{\sc H.~Wilson and J.~Cowan}, {\em Excitatory and inhibitory interactions in
  localized populations of model neurons}, Biophys. J., 12 (1972), pp.~1--24.

\bibitem{wilson-cowan:73}
\leavevmode\vrule height 2pt depth -1.6pt width 23pt, {\em A mathematical
  theory of the functional dynamics of cortical and thalamic nervous tissue},
  Biological Cybernetics, 13 (1973), pp.~55--80.

\bibitem{woods-champneys:99}
{\sc P.~Woods and A.~Champneys}, {\em Heteroclinic tangles and homoclinic
  snaking in the unfolding of a degenerate reversible hamiltonian-hopf
  bifurcation}, Physica D: Nonlinear Phenomena, 129 (1999), pp.~147--170.

\end{thebibliography}

\end{document}